\documentclass[11pt]{article}

\usepackage[margin=1.1in]{geometry}
\usepackage{hyperref}
\usepackage{cite}
\usepackage{color}

\definecolor{darkpastelgreen}{rgb}{0.2, 0.55, 0.24}
	\definecolor{cadmiumgreen}{rgb}{0.0, 0.42, 0.24}
\definecolor{armygreen}{rgb}{0.29, 0.33, 0.13}
\hypersetup{colorlinks,linkcolor={blue},citecolor={darkpastelgreen},urlcolor={red}}

\usepackage{cite}
\usepackage[round]{natbib}

\usepackage[utf8]{inputenc} % allow utf-8 input
\usepackage[T1]{fontenc}    % use 8-bit T1 fonts
\usepackage{hyperref}       % hyperlinks
\usepackage{url}            % simple URL typesetting
\usepackage{booktabs} 
\usepackage{tabularx}
\usepackage{wrapfig}

\usepackage{amsfonts}       % blackboard math symbols
\usepackage{nicefrac}       % compact symbols for 1/2, etc.
\usepackage{microtype}      % microtypography
\usepackage[dvipsnames]{xcolor}
\usepackage{microtype}
\usepackage{graphicx}
\usepackage{subfigure}
\usepackage{booktabs} % for professional tables
\usepackage{algorithmicx}
\usepackage{algorithm}% http://ctan.org/pkg/algorithms
\usepackage{algpseudocode}% http://ctan.org/pkg/algorithmicx
\usepackage{xtab,booktabs}
\usepackage{array, ltablex, multirow}
\usepackage{amsmath,amssymb,amsthm}
\usepackage{bbm}
\usepackage{xcolor}
\usepackage{graphicx}
\usepackage{subfigure}
\usepackage[utf8]{inputenc} % allow utf-8 input
\usepackage[T1]{fontenc}    % use 8-bit T1 fonts
\usepackage{hyperref}       % hyperlinks
\usepackage{url}            % simple URL typesetting
\usepackage{booktabs}       % professional-quality tables
\usepackage{amsfonts}       % blackboard math symbols
\usepackage{nicefrac}       % compact symbols for 1/2, etc.
\usepackage{amssymb}
\usepackage{enumitem}

\usepackage{pifont}% http://ctan.org/pkg/pifont
\DeclareMathOperator*{\argmax}{arg\,max}
\DeclareMathOperator*{\argmin}{arg\,min}

\let\proof\relax \let\endproof\relax

\title{Minimax Optimization:  The Case of Convex-Submodular}

\author{Arman Adibi\thanks{Department of Electrical and Systems Engineering, University of Pennsylvania, Philadelphia, PA, USA. } ,\quad  Aryan Mokhtari\thanks{Department of Electrical and Computer Engineering, The University of Texas at Austin, Austin, TX, USA.} ,\quad Hamed Hassani\footnotemark[1]}

\newcommand{\0}{\mathbf{0}}

\newcommand{\K}{\mathcal{K}}

\newcommand{\X}{\mathcal{X}}

\newcommand{\I}{\mathcal{I}}

\newcommand{\x}{\mathbf{x}}
\newcommand{\y}{\mathbf{y}}

\newcommand{\cmark}{\ding{51}}%
\newcommand{\xmark}{\ding{55}}%
\newtheorem{theorem}{Theorem}
\newtheorem{corollary}{Corollary}
\newtheorem{lemma}{Lemma}
\newtheorem{proposition}{Proposition}
\newtheorem{assumption}{Assumption}

\newtheorem{defi}{Definition}

\newtheorem{remark}{Remark}

\begin{document}


\maketitle

\begin{abstract}
Minimax optimization has been central in addressing various applications in machine learning, game theory, and control theory. Prior literature has thus far mainly focused on studying such problems in the continuous domain, e.g., convex-concave minimax optimization is now understood to a significant extent.  Nevertheless, minimax problems extend far beyond the continuous domain to mixed continuous-discrete domains or even fully discrete domains. In this paper, we study mixed continuous-discrete minimax problems where the minimization is over a continuous variable belonging to Euclidean space and the maximization is over subsets of a given ground set. 
We introduce the class of convex-submodular minimax problems, where the objective is convex with respect to the continuous variable and submodular with respect to the discrete variable.  Even though such problems appear frequently in machine learning applications, little is known about how to address them from algorithmic and theoretical perspectives. 
For such problems, we first show that obtaining saddle points are hard up to any approximation, and thus introduce new notions of (near-) optimality. We then provide several algorithmic procedures for solving convex and monotone-submodular minimax problems and characterize their convergence rates, computational complexity, and quality of the final solution according to our notions of optimally. Our proposed algorithms are iterative and combine tools from both discrete and continuous optimization. 
Finally, we provide numerical experiments to showcase the effectiveness of our purposed methods.
\end{abstract}

%%%%%%%%%%%%%%%%%%%%%%%%%%%%%
\section{Introduction}

The problem of solving a minimax optimization problem, also known as the saddle point problem, appears in many domains such as robust optimization \citep{ben2009robust}, game theory \citep{osborne1994course}, and robust control~\citep{zhou1998essentials,hast2013pid}. It has also recently attracted a lot of attention in the machine learning community due to the rise of  generative adversarial networks (GANs) \citep{goodfellow2014generative} and robust learning \citep{bertsimas2011theory,lanckriet2002robust,li2019robust}.
There has been an extensive literature on the design of convergent methods for solving minimax problems for the case that both minimization and maximization variables belong to continuous domains \citep{tseng1995linear,nesterov2007dual,li2015accelerated,ouyang2019lower,thekumparampil2019efficient,zhao2019optimal,hamedani2018primal,alkousa2019accelerated,daskalakis2017training,ibrahim2020linear,nouiehed2019solving,mokhtari2020unified,lin2020gradient,murty1985some}. In particular, for the case that the loss function is (strongly) convex with respect to the minimization variable and (strongly) concave with respect to maximization variable several efficient algorithms have been studied \citep{nesterov2007dual,li2015accelerated,mokhtari2020convergence}, including the extra-gradient method \citep{korpelevich1976extragradient,nemirovski2004prox} that is known to be optimal for this setting.
Subsequently, multiple other algorithms, such as Nesterov’s dual extrapolation \citep{nesterov2007dual}, accelerated proximal gradient \citep{li2015accelerated}, and optimistic gradient descent-ascent \citep{mokhtari2020convergence} methods have been introduced and analyzed for solving minimax problems.
However, all these methods suffer from two major limitations: (i) they are provably convergent only in convex-concave settings; (ii) they are designed for the settings that both minimization and maximization variables belong to continuous domains. 

There has been some effort to address the first limitation by finding a first-order stationary point or locally stable point for the problems that are not convex-concave \citep{lin2020near,diakonikolas2021efficient,yang2020global,sanjabi2018solving}. However, these approaches fail to guarantee any global optimality as it is known that finding a saddle point in a general nonconvex-nonconcave setting is NP-hard \citep{jin2020local}. Nonetheless,  it might be possible to achieve global approximation guarantees for \emph{structured} saddle minimax problems. Addressing the second limitations and developing methods for discrete-continuous domains or fully discrete domains requires exploiting tools from discrete optimization. Several recent works have considered applications involving specific discrete-continuous minimax problems and proposed structure-informed algorithms \citep{zhou2018minimax}. However, to our knowledge, there is no work that provides a principled algorithmic or theoretical framework to study minimax problems with mixed discrete-continuous components and it is not even clear if such problems allow for tractable solutions with global guarantees.  
\begin{table*}[t!]
\caption{Algorithms performance guarantee. Here $c_{f}$ is the cost of single computation of $f$, $c_{P_{\x}}$ and $c_{P_{\y}}$ are cost of projection in $\mathcal{X}$ and $\mathcal{Y}$, $c_{{\nabla_\x} f }$ is the cost of computing gradient of $f$ with respect to $\x$, and $c_{{\nabla_\x} F }$ and $c_{{\nabla_\y} F }$
 are the cost of computing gradient of multilinear extension $F$ with respect to $\x$  and $\y$, respectively. $k$ is the cardinality constraint ($|S|\leq k$) and $n$ is size of the ground set $|V|=n$. \label{tab:1}
 }
 \vspace{2mm}
 
\centering

\begin{tabular}{|p{0.95cm}|p{1.2cm}|p{1.1cm}|p{4.5cm}|p{1.1cm}|p{1.2cm}|p{1.7cm}|} 
  \hline
 Alg. & Number of iterations & Approx.\quad ratio & Cost per iteration& Card. const. & Matroid const.  &  Unbounded grad.\\
 \hline
 GG   & $\mathcal{O}({1}/{\epsilon^2})$  & $1-1/e$ & \small{$nk.c_{f}+c_{P_{\x}}+c_{{\nabla_\x} f  }$}&\cmark&\xmark & \xmark \\
 \hline
 GG   & $\mathcal{O}({1}/{\epsilon^2})$  & $1/2$ & $nk.c_{f}+c_{P_{\x}}+c_{{\nabla_\x} f }$&\xmark&\cmark & \xmark \\
\hline
 GRG    & $\mathcal{O}({1}/{\epsilon^2})$ & $1/2$ & $(n+k)c_{f}+c_{{\nabla_\x} f +c_{P_{\x}}}$&\cmark&\xmark & \xmark \\
\hline
EGG   &$\mathcal{O}({1}/{\epsilon^2})$  & $1-1/e$ & $2nk.c_{f}+2c_{P_{\x}}+2c_{{\nabla_\x} f }$&\cmark&\xmark & \cmark   \\
\hline
EGG   &$\mathcal{O}({1}/{\epsilon^2})$  & $1/2$ & $2nk.c_{f}+2c_{P_{\x}}+2c_{{\nabla_\x} f }$&\cmark&\cmark & \cmark   \\
\hline
 EGRG & $\mathcal{O}({1}/{\epsilon^2})$ & ${1}/{2}$& $2(n+k)c_{f}+2c_{P_{\x}}+2c_{{\nabla_\x} f }$ &\cmark&\xmark & \xmark\\
\hline
  EGCE &$\mathcal{O}({1}/{\epsilon})$ & $1/2$& $2c_{P_{\x}}+2c_{P_{\y}}+2c_{{\nabla_\x} F }+2c_{{\nabla_\y} F }$ &\cmark&\cmark & \cmark  \\
 \hline
\end{tabular}
\end{table*}

%%%%%%%%%%%%%%%%%%%%%%%%%%%%%%%%%%%%%%%%%%%%%%%%%%%%%%%%%%%%%%%%%%%%%%%%%%%%%%%%%%%%%%%%%%%%%%%%%%%%%%%%%%%%%%%%%%%%%%%%%%%%%%%%%%%%%%%%%%%%%%%%%%%%%%%%%%%%%%%%%%%%%%%%%%

In this paper, we tackle these two issues and present iterative methods with theoretical guarantees to solve structured non convex-concave minimax problems, where the minimization variable is from a continuous domain and the maximization variable belongs to a discrete domain. Concretely, for a non-negative function $f:\mathbb{R}^d \times 2^V \to \mathbb{R}_{+}$, consider the minimax problem 
%%%%
\begin{equation}\label{eq:main_problem}
   {\rm{OPT}} \triangleq \min_{\x \in \mathcal{X}}\max_{S \in \mathcal{I}} f(\x,S),
\end{equation}
where $\x$ belongs to a convex set $\mathcal{X}\subset \mathbb{R}^d$ and $S$  is a subset of the ground set ${V}$ with $n$ elements that is constrained to be inside a matroid~$\mathcal{I}$. Given a fixed $S$, the function $f(\cdot,S)$ is convex with respect to the continuous (minimization) variable. Further, given a fixed $\x$, the  function $f(\x,\cdot)$ is submodular with respect to the discrete (maximization) variable. We refer to this problem as \emph{convex-submodular minimax problem}.
%The loss function $f$ is convex with respect to $\x$ and submodular with respect to $S$.

The convex-submodular minimax problem in \eqref{eq:main_problem} encompasses various applications. In Section~\ref{sec:Poisoning Attack}, we describe specific optimization problems, such as convex-facility-location, as well as applications such as designing adversarial attacks on recommender systems. There are various other applications that can be cast into Problem~\eqref{eq:main_problem}, in particular, when convex models have to be learned while data points are selected or changed according to notions of summarization, diversity, and deletion. Examples include  learning under data deletion~\citep{ginart2019making, neel2020descent, wu2020deltagrad}, robust text classification~\citep{lei2018discrete}, minimax curriculum learning~\citep{zhou2018minimax,zhou2021curriculum,zhou2020curriculum,soviany2021curriculum}, minimax supervised learning~\citep{farnia2016minimax}, and minimax  active learning~\citep{ebrahimi2020minimax}. \

\subsection{Our Contributions}\label{sec:contributions}

In this paper, we provide a principled study of the problem defined in \eqref{eq:main_problem}, from both theoretical and algorithmic perspectives, when $f$ is convex in the minimization variable and submodular as well as \emph{ monotone}  in the maximization variable\footnote{For completeness, a function $g : 2^V \to \mathbb{R}$ is called submodular if for any two subsets $S, T \subseteq V$ we have: \; $g(S\cap T) + g(S\cup T) \leq g(S) + g(T)$.  Moreover, $g$ is called monotone if for any $S\subseteq T$ we have $g(S)\leq g(T)$.}.
We  introduce efficient iterative algorithms for solving this problem and develop a theoretical framework for analyzing such algorithms with guarantees on the quality of the resulting solutions according to the notions of optimality that we define.

\vspace{3mm}
\noindent\textbf{Notions of (near-)optimality and hardness results.} For minimax problems, the strongest notion of optimality is defined through saddle points or their approximate versions. We first provide a negative result that shows finding a saddle point or any approximate version of it (which we term as an $(\alpha,\epsilon)$-saddle point) is NP-hard for general convex-submodular problems (Theorem~\ref{thm:saddle_neg}).
We thus introduce a slightly weaker notion of optimality that we call $(\alpha,\epsilon)$-approximate minimax solutions for Problem~\eqref{eq:main_problem}. Roughly speaking, the quality of the minimax objective at such solutions is at most $\frac{1}{\alpha}(\text{OPT} + \epsilon)$, and hence they are near-optimal when $\alpha < 1$. We show in Theorem~\ref{thm:approx_neg} that obtaining such solutions for $\alpha > 1-1/e$ is NP-hard. This is a non-trivial result that does not readily follow from known hardness results in submodular maximization. Consequently, we focus on efficiently finding solutions in the regime of $\alpha\leq (1-1/e)$. We present several algorithms that achieve this goal and theoretically analyze their complexity and quality of their solution. 

\vspace{3mm}
\noindent\textbf{Algorithms with guarantees on convergence rate,  complexity, and solution quality.} Our proposed algorithms are as follows (see also Table~\ref{tab:1}): (i) \emph{Greedy-based methods.} We first present Gradient-Greedy (GG), a method alternating between gradient descent for minimization and greedy for maximization. We further introduce Extra-Gradient-Greedy (EGG) that uses an extra-gradient step instead of gradient step for the minimization variable. We prove that both algorithms achieve a
$((1-1/e),\epsilon)$-approximate minimax solution after $\mathcal{O}(1/\epsilon^2)$ iterations when $\mathcal{I}$ is a cardinality constraint. Importantly, EGG does not require the bounded gradient norm condition as opposed to GG. Our results for the case that $\mathcal{I}$ is a matroid constraint (see Table~\ref{tab:1}) are provided in the Appendix. (ii)~\emph{ Replacement greedy-based methods.} The greedy-based methods require $\mathcal{O}(nk)$ function computations at each iteration. To improve this complexity, we  present alternating methods that use replacement greedy for the maximization part to reduce the cost of each iteration to $\mathcal{O}(n)$.  {The Gradient Replacement-Greedy (GRG) algorithm achieves a $(1/2,\epsilon)$-approximate minimax solution after $\mathcal{O}(1/\epsilon^2)$ iterations and Extra-Gradient Replacement-Greedy (EGRG) achieves a $({1}/{2},\epsilon)$-approximate minimax solution after $\mathcal{O}(1/\epsilon^2)$, when $\mathcal{I}$ is a cardinality constraint.} (iii)~\emph{Continuous extension-based methods.} Note that all mentioned methods  achieve a convergence rate of  $\mathcal{O}(1/\epsilon^2)$. To improve this convergence rate, we further introduce the extra-gradient on continuous extension (EGCE) method that runs extra-gradient update on the continuous extension of the submodular function. We show that EGCE is able to achieve an  $({1}/{2},\epsilon)$-approximate minimax solution after at most $\mathcal{O}(1/\epsilon)$ iterations, when $\mathcal{I}$ is a general matroid constraint.

\subsection{Related Work}

%{\color{red}{we need to rewrite the related work section}}
%\iffalse
%\textbf{Submodular Maximization.}

%Many applications in machine learning such as recommendation systems\citep{krause2014submodular,mitrovic2018data,adibi2020submodular}, data summerization\citep{kirchhoff2014submodularity,wei2013using}, and active learning\citep{golovin2011adaptive} can be formulated as submodular maximization. This problem is NP-hard in general but you can achieve $1-1/e$ approximate solution with Greedy algorithm\citep{krause2014submodular, nemhauser1978best,wolsey1982analysis}. Many other extension of this problem has been proposed later on such as online submodular maximization\citep{jegelka2011online,streeter2009online,golovin2014online,chen2018projection}, distributed submodular maximization\citep{mirzasoleiman2013distributed,mirzasoleiman2016distributed}, robust submodular maximization\citep{krause2008robust,iyer2021unified,orlin2018robust}, stochastic submodular optimization \citep{karimi2017stochastic}, K-submodular optimization\citep{ohsaka2015monotone}, two-stage submodular optimization\citep{mitrovic2018data}, decentralized submodular optimization\citep{mokhtari2018decentralized}, deletion robust submodular optimization\citep{mirzasoleiman2017deletion}, and multi-agent Submodular Optimization\citep{santiago2016multivariate}. \fi 
Several recent works have considered specific applications that require solving Problem~\eqref{eq:main_problem} when $f$ is nonconvex-submodular \citep{zhou2018minimax,lei2018discrete,mirzasoleiman2020coresets}. \cite{zhou2018minimax} consider the problem of minimax curriculum learning and propose an algorithm similar to gradient-greedy (GG). They provide an upper bound on the distance between their obtained solution and the optimal solution when $f$ is strongly convex in $\x$ and monotone-submodular in $S$ with non-zero curvature.  Moreover, \cite{lei2018discrete} study designing an adversarial attack in text classification and show that for some specific neural network structures, the task of designing an adversarial attack can be formulated as submodular maximization, leading to a minimax nonconvex-submodular problem. An algorithm similar to gradient-greedy is then proposed by \cite{lei2018discrete} for designing attacks and it has led to  successful experimental results. In contrast, this paper is the first to introduce a principled study of Problem~\eqref{eq:main_problem} for general functions $f$ with newly developed notions of optimality, algorithmic frameworks, and theoretical guarantees. 
 
 %Also, minimax Curriculum Learning\citep{zhou2018minimax} is minimax noncovex-submodular problem. They consider the problem of Curriculum learning, which is a problem of designing a sequence of training sets that can improve the learning algorithm. For this specific application they proposed similar algorithm to gradient-greedy and analyzed it for the strongly convex-submodular case. Non of the above mentioned works proposed,analyzed and studied the convergence of algorithms we have for minimax convex-submodular problem.

Another relevant line of work is the literature on robust submodular optimization \citep{krause2008robust,bogunovic2017robust,mirzasoleiman2017deletion,kazemi2018scalable,bogunovic2018robust,iyer2021unified,orlin2018robust,chen2017robust,anari2019structured,wilder2018equilibrium,bogunovic2017distributed,mitrovic2017streaming}. This setting corresponds to solving a \emph{max-min} optimization problem which involves only \textit{discrete variables}, and hence, it is  different from our setting with fundamentally different methods. 
For such problems, finding discrete solutions with any approximation factor is NP-hard; and consequently, the literature has mostly focused on obtaining solutions that satisfy a bi-criteria approximation guarantee. 
%The max-min version of Problem~\eqref{eq:main_problem} can be related to this line of work in some special cases and through appropriate continuous relaxations of the discrete variables. 
Another related work is distributionally robust submodular maximization in \citep{staib2019distributionally} which  is  a special case of \emph{max-min} version of Problem~\eqref{eq:main_problem}. In this setting, the inner minimization has a special structure that allows for a closed form solution, and hence, the problem can be solved by using appropriate techniques from continuous submodular optimization. We will derive the implication of our results on the max-min version of Problem~\eqref{eq:main_problem} in the Appnedices. 
%Again, this setting and the corresponding techniques/guarantees do not
%We emphasize here that in this paper wefocus mainly on the Problem~\eqref{eq:main_problem}, and since for this problem obtaining (approximate) saddle-point is hard, the tools that we develop are fundamentally different from the ones developed in prior work.

%%%%%%%%%%%%%%%%%%%%%%%%%%%%

\section{Convex-Submodular Minimax Optimization}

For the minimax problem in~\eqref{eq:main_problem}, a natural goal is to find a so-called saddle point.
%which can be interpreted as an equilibrium if we consider the minimization and maximization problems as a game between two players. 
Next, we formally define the notion of saddle point for Problem~\eqref{eq:main_problem}.

%%%%%%%%%%%%%%%%%%%%%%%%%%%%
\begin{defi}
A pair $(\x^*, S^*)$ is a saddle point of the function $f$ if the following condition holds:
%%%%%%%%%%%%%%%%%%%%%%%%%%%%
\begin{equation}\label{saddle_def}
    \forall \x \in \mathcal{X}, S \in \mathcal{I}: \,  f(\x^*, S) \leq f(\x^*, S^*) \leq f(\x,S^*)
\end{equation}
%%%%%%%%%%%%%%%%%%%%%%%%%%%%
\end{defi}
%%%%%%%%%%%%%%%%%%%%%%%%%%%%
\vspace{2mm}

Based on this definition, $(\x^*, S^*)$ is a saddle point of Problem \eqref{eq:main_problem}, if there is no incentive to modify the minimization variable $\x^*$ when the maximization variable  is fixed and equal to $S^*$, and, conversely, there is no incentive to change the maximization variable from $S^*$ when the minimization variable is $\x^*$. 
In other words, a saddle point can be interpreted as an equilibrium.

There is a rich literature on efficient approaches for finding an $\epsilon$-saddle point for convex-concave minimax optimization, where $\epsilon$ is an arbitrary positive constant \citep{thekumparampil2019efficient}. To define an $\epsilon$-saddle point, we first need to define the duality gap, which is given by
$
D(\x,S) := \bar{\phi}(\x) - \underline{\phi}(S)$,
where 
$$
\bar{\phi}(\x)\!\!:=\max_{S\in \I}f(\x,S), \qquad \underline{\phi}(S):=\min_{\x\in \mathcal{X}}f(\x,S).
$$
Considering these definitions, we call a pair of solution $\epsilon$-saddle point if their duality gap is at most $\epsilon$.

{\defi
A pair $(\bar\x,\bar S)$ is called an  $\epsilon-$saddle point of $f$ if it satisfies
%%%%%%%%%%%%%%%%%%%%%%%%%%%%
\begin{equation}\label{eps_saddle}
   D(\bar\x,\bar S)= \bar{\phi}(\bar\x)-\underline{\phi}(\bar S)\leq \epsilon.
\end{equation}
}
\vspace{2mm}

One can verify that if we set $\epsilon=0$, then Definitions~1 and 2 coincide, i.e., $(\x^{*},S^{*})$ satisfies \eqref{eps_saddle} for $\epsilon=0$ if and only if $(\x^{*},S^{*})$ satisfies \eqref{saddle_def}. Hence, to derive a finite time analysis we often aim for finding an $\epsilon$-saddle point. For instance,  for smooth and convex-concave problems extra-gradient obtains an $\epsilon$-saddle point after  $\mathcal{O}(1/\epsilon)$ iterations (which is the optimal complexity). 
%which is the optimal complexity for smooth convex-concave problems as follows.

However, for our convex-submodular setting, one cannot expect to find an $\epsilon$-saddle point efficiently, as the special case of finding an $\epsilon$-accurate solution for submodular maximization is in general NP-hard \citep{ nemhauser1978best,wolsey1982analysis,krause2014submodular}. Although solving the problem of maximizing a  monotone submodular function subject to a matroid constraint is hard, one can find $\alpha$-approximate solution of that in polynomial time, i.e., finding a solution that its function value is at least $\alpha \rm{OPT}$, where $\alpha \in (0,1)$. Inspired by this observation, we introduce the notion of $(\alpha,\epsilon)$-saddle point for our convex-submodular setting.

%%%%%%%%%%%%%%%%%%%%%%%%%%%%
%\begin{lemma}Finding a Saddle point of problems \eqref{eq:main_problem} and \eqref{eq:main_problem2} is NP-hard.
%\end{lemma}
%%%%%%%%%%%%%%%%%%%%%%%%%%%%

%%%%%%%%%%%%%%%%%%%%%%%%%%%%
%{\defi 
%Define $\bar{\phi}(\x)\!\!:=\max_{S\in \I}f(\x,S)$ and $\underline{\phi}(S):=\min_{\x\in \mathcal{X}}f(\x,S)$.}
%%%%%%%%%%%%%%%%%%%%%%%%%%%%
{\defi
A pair $(\hat{\x},\hat{S})$ is called an  $(\alpha,\epsilon)-$saddle point of $f$ if it satisfies
%%%%%%%%%%%%%%%%%%%%%%%%%%%%
\begin{equation}\label{eq:saddle_def}
    \alpha\bar{\phi}(\hat{\x})-\underline{\phi}(\hat{S})\leq \epsilon.
\end{equation}
\label{def:saddle}}
\vspace{2mm}

Our first result is a negative result that shows even finding an $(\alpha,\epsilon)$-saddle point is not tractable. 

\begin{theorem}\label{thm:saddle_neg}
Finding $(\alpha,\epsilon)$-saddle point for Problem~\eqref{eq:main_problem} is NP-hard for any  $\alpha>0$.
\end{theorem}
%%%
While this result shows intractability of finding (approximate) saddle-points for Problem~\eqref{eq:main_problem}, one avenue to provide solutions with guaranteed quality is to see whether we can find solutions that achieve a fraction of OPT. We thus proceed to introduce the notion of \emph{approximate minimax solution}.

{\defi
We call a point $\hat{\x}$ an $(\alpha,\epsilon)$-approximate minimax solution of Problem~\eqref{eq:main_problem} if it satisfies
%%%%%%%%%%%%%%%%%%%%%%%%%%%%
\begin{equation}\label{eq:(alpha,eps)sol}
  \alpha\bar{\phi}(\hat{\x})\leq {\rm{OPT}}+{\epsilon},
\end{equation}
where ${\rm{OPT}}= \min_{\x\in \mathcal{X}}\bar{\phi}(\x) =\min_{\x \in \mathcal{X}}\max_{S \in \mathcal{I}} f(\x,S)$.
\label{def:approx_sol}}
\vspace{2mm}

Next, we describe the notion of an  $(\alpha,\epsilon)$-approximate minimax solution for Problem~\eqref{eq:main_problem}.
The minimax problem in \eqref{eq:main_problem} can be interpreted as a sequential game, where we first select an action $\x$ and then an adversary chooses a set $S$ to maximize our loss $f(\x,S)$. In this case, our goal is to find $\x$ that minimizes the loss obtained by the worst possible action by the adversary, i.e., we aim to minimize the function $\bar{\phi}(\hat{\x}):=\max_{S\in \I}f(\hat{\x},S)$  over the choice of $\hat{\x}$. Indeed, finding the exact minimizer is also hard and we should seek approximate solutions. Hence, our goal is to find solutions $\hat{\x}$ whose worst-case loss $\bar{\phi}(\hat{\x})$ is only a factor larger than the best possible loss $\text{OPT}=\min_{\x \in \mathcal{X}} \bar{\phi}({\x})$. That said, by finding an  $(\alpha,\epsilon)$-approximate minimax solution for Problem~\eqref{eq:main_problem} we obtain a solution whose loss is at most $ (\text{OPT}+\epsilon)/{\alpha}$, where $0\!<\!\alpha\!\leq\! 1$ and $\epsilon>0$.

The task of finding an $\hat{\x}$ that is  $(\alpha,\epsilon)$-approximate minimax solution is easier than finding a pair $(\tilde{\x},\tilde{S})$ that is an $(\alpha,\epsilon)$-saddle point, since if the pair $(\tilde{\x},\tilde{S})$ satisfies  \eqref{eq:saddle_def}, then $\tilde{\x}$ satisfies    \eqref{eq:(alpha,eps)sol}:
\begin{align*}
      &\alpha\bar{\phi}(\tilde{\x})-\underline{\phi}(\tilde{S})\leq \epsilon
      \quad
      \Rightarrow
      \quad
       \alpha\bar{\phi}(\tilde{\x})-\max_{S\in \mathcal{I}}\underline{\phi}(S)\leq \epsilon
        \quad 
     \Rightarrow
      \quad
       \alpha\bar{\phi}(\tilde{\x})-\min_{\x\in \mathcal{X}}\bar{\phi}(\x)\leq \epsilon
\end{align*}
Hence, the condition in \eqref{eq:saddle_def} is more strict compared to \eqref{eq:(alpha,eps)sol}. In fact, in the next section, we show that unlike the task of finding an $(\alpha,\epsilon)$-saddle point of Problem~\eqref{eq:main_problem} that is NP-hard for any $\alpha\in (0,1]$, one can find an $(\alpha,\epsilon)$-approximate minimax solution of Problem~\eqref{eq:main_problem} in poly-time for $\alpha\in (0,1-1/e]$. Alas, the problem is still NP-hard for $\alpha\in (1-1/e,1]$ as we show in Theorem~\ref{thm:approx_neg}. 

%%%%%%%%%%%%%%%%%%%%%%%%%%%%
\begin{theorem}\label{thm:approx_neg}
 Let $\alpha=1-1/e+\gamma$ for a positive constant $\gamma>0$. If there exists a polynomial time algorithm and a polynomial time oracle that can achieve an $(\alpha,\epsilon)$-approximate solution for any choice of the function $f(\x,S)$ in problem \eqref{eq:main_problem}, then P = NP.
\end{theorem}

We emphasize that Theorem~\ref{thm:approx_neg} does not follow directly from that fact that submodular maximization beyond $(1-1/e)$-approximation is hard, and hence it is non-trivial. Indeed, one naive way to argue for the proof of this theorem (which is incorrect) is to consider functions $f(\x,S)$ whose output does not depend on the variable $\x$, i.e. $f(\x,S) = f(S)$, and  use the hardness results for submodular optimization. But for such functions \emph{any} point $\x$ is an optimal solution (with $\alpha=1$). Hence, the proof of the theorem (provided in the appendix) requires a novel idea beyond trivial consequences of known results for submodularity.

So far we have shown two results: (i) Finding an approximate $(\alpha, \epsilon)$-saddle point is hard for $\alpha\!>\!0$. (ii)~We introduced the notion of $(\alpha, \epsilon)$-approximate solution and showed that for $\alpha\!>\! 1\!-\!1/e$ finding an $(\alpha, \epsilon)$-approximate solution is hard. 
The only missing piece is showing whether or not it is possible to \emph{efficiently} find an $(\alpha,\epsilon)$-approximate minimax solution when $\alpha\!\leq\! 1\!-\!1/e$. In the rest of the paper, we provide an affirmative answer to this question and present methods achieving this goal.

%%%%%%%%%%%%%%%%%%%%%%%%%
\section{Algorithms}\label{sec:algs}

In this section, we present a set of algorithms that are able to find an $(\alpha,\epsilon)$-approximate minimax solution of Problem~\eqref{eq:main_problem}. To present these algorithms, we first present two subroutines that we use in the implementation of our algorithms\footnote{For better exposition,  we consider the case that $\mathcal{I}$ is   $k$-carnality constraint and refer to  Appendix for~matroids.}: (i)
 greedy update and (ii) replacement greedy update. 
 
    \vspace{2mm}
\noindent
\textbf{Greedy subroutine.} In the greedy update, for a fixed minimization variable $\x$, we select a subset $S$ with $k$ elements in a greedy fashion, i.e., we sequentially pick $k$ elements that maximize the marginal gain. Specifically, if we define $\Delta_e f( \x,S)=f( \x,S\cup \{e\})-f( \x,S)$ as the marginal gain of element $e$, in the greedy update, for a given variable $\x$ we perform the update
\begin{align}\label{greedy_update}
 S_{i+1}=S_{i}\cup  \{\arg \max\Delta_e f( \x,S)\}   ,
\end{align}
for $i=0,\dots, k-1$, where $S_0$ is the empty set. The output of this process is $S_k$ with $k$ elements. We use the notation \textsc{Greedy}$(f,k,\x)$ for the greedy subroutine, which takes function $f$, cardinality constraint parameter $k$, and variable $\x$ as inputs, and returns a set $S$ by performing \eqref{greedy_update} for $k$ steps.

    \vspace{2mm}
\noindent \textbf{Replacement greedy subroutine.} In the replacement greedy update \citep{mitrovic2018data,schrijver2003combinatorial,pmlr-v70-stan17a}, for a given variable $\x$ and set $S$, the output is an updated set $S^+$ whose function value at $\x$ is larger than the one for  $S$, i.e., $f(\x,S)\leq f(\x,S^+)$. The procedure for finding the new set $S^+$ is relatively simple. If the size of the input set $S$ is less than $k$, we add one more element to the set $S$ that maximizes the marginal gain and the resulted set would be $S^+$. In other words, if $|S|<k$, 
\begin{equation}\label{rep_case_1}
 S^+=S\cup  \{\arg \max\Delta_e f( \x,S)\}.
 \end{equation}
If the size of the input set $S$ is $k$, we first remove one element of the set $S$ that leads to minimum decrease in the function value (denoted by $e^*$) and then replace it with another element of the ground set that maximizes the marginal gain. Hence, if $|S|=k$, we have
\begin{equation}\label{rep_case_2}
 S^+=(S\setminus \{e^*\}) \cup  \{\arg \max\Delta_e f( \x,S\setminus \{e^*\})\},
   \end{equation}
   where
$ e^*=\argmax_{e\in  S}\{ f( \x ,S\setminus e)\}$.

We use the notation \textsc{RepGreedy}$(f,k,\x,S)$ for the replacement greedy subroutine. 
%which takes function $f$, cardinality constraint parameter $k$, minimization variable $\x$, and set $S$ as inputs, and returns a set $S^+$ that is computed according to \eqref{rep_case_1} or \eqref{rep_case_2}, depending on the size of the input set $S$. 
Note that replacement greedy is computationally cheaper than greedy, as it requires only one pass over the ground set, while greedy requires $k$ passes.
%%%%%%%%%%%%%%%%%%%%%alg%%%%%%%%%%%%%%%%%%%%%%%%
%\begin{algorithm}[H]
%\caption{\textsc{GREEDY}$(f,\mathcal{I},\x)$}\label{alg:greedy}
%\begin{algorithmic}
%\State \textbf{Initialize} The set $S_0$ to the empty set and let $k$ to be the rank of matroid $\I$.
%\For{$i=0$ to $k-1$} 
%%%%%%%%%%%%%%%%%%%%%%%
%\vspace{-3mm}
% \State $$S_{i+1}=S_{i}\cup  \{\arg \max\Delta_e f( \x,S):S_{i}\cup\{e\}\in \mathcal{I}\}$$
% \vspace{-4mm}
% \EndFor
% \State return $S_k$
%\end{algorithmic}%
%\end{algorithm}
%%%%%%%%%%%%%%%%%%%%%%%ALG%%%%%%%%%%%%%%%%%%%%%%%
%%%%%%%%%%%%%%%%%%%%%%%%%%%%%%%%%%%%%%%%%%%%%%%%%
\subsection{Greedy-based Algorithms}
Next, we present greedy-based methods to find $(\alpha, \epsilon)$-approximate minimax solutions for Problem~\eqref{eq:main_problem}.

    \vspace{2mm}
\noindent{\textbf{Gradient Greedy.}} The first algorithm that we present is Gradient Greedy (GG), which uses a projected gradient descent step to update the minimization iterate $\x_t$ at each iteration, i.e., 
$
{\x}_{t+1} =  \pi_{\mathcal{X}}(\x_{t} -\gamma_t  \nabla f(\x_{t},S_t)),
$
and then uses a greedy procedure to update the maximization variable $S_t$. This update is performed in an alternating fashion, where we first use $\x_t$ and $S_t$ to find $\x_{t+1}$ and then we use the updated variable $\x_{t+1}$ to compute $S_{t+1}$. Note that the final output of this process is a weighted average of all variables $\x_t$ that are observed from time $t=1$ to $t=T$, defined as ${\x}_{sol}=({\sum_{t=1}^{T}\gamma_t})^{-1}{\sum_{t=1}^{T}\gamma_t\x_t}$. The steps of GG are summarized in Algorithm~\ref{alg:(Extra)-gradient-greedy} {option~\text{I}}.

%%%%%%%%%%%%%%%%%%%%%alg%%%%%%%%%%%%%%%%%%%%%%%%
\begin{algorithm}[t]
\caption{\\\text{\textbf{Option \text{I}:} Gradient Greedy (GG)}\quad \\\textbf{Option \text{II}:} Extra-Gradient Greedy (EGG)}\label{alg:(Extra)-gradient-greedy}

\begin{algorithmic}
\State \textbf{Initialize} the set $S_1$ to $\emptyset$ and variable $\x_1$ to zero. 
\For{$t=1$ to $T$} 
%%%%%%%%%%%%%%%%%%%%%%%
%%%%%%%%%%%%%%%%%%%%%%%
\State$\textbf{Option  \text{I}:}$ 
   $
{\x}_{t+1} =  \pi_{\X}(\x_{t} -\gamma_t  \nabla f(\x_{t},S_t))$
\State \qquad\qquad\,\,\,\,\,\,$
{S}_{t+1} = \textsc{Greedy}(f , k,{\x}_{t+1})
$
\State$\textbf{Option  \text{II}:} $$\ 
\hat{\x}_t =  \pi_{\X}(\x_t -\gamma_t  \nabla f(\x_t,S_t))
$

\State \qquad\qquad\,\,\,\,\,\,$
\hat{S}_t = \textsc{Greedy}(f , k,\hat{\x}_t)
$

\State \qquad\qquad\,\,\,\,\,\,$
\x_{t+1} = \pi_{\X}( \x_t -\gamma_t  \nabla f(\hat{\x}_t,\hat{S}_t))
$
\State \qquad\qquad\,\,\,\,\,\,$
{S}_{t+1} = \textsc{Greedy}(f , k,{\x}_{t+1})
$

\EndFor
%%%%%%%%%%%%%%%%%%%%%%%
\State $\textbf{Option  \text{I}:}$ ${\x}_{sol}=({\sum_{t=1}^{T}\gamma_t})^{-1}{\sum_{t=1}^{T}\gamma_t\x_t}$
\State$\textbf{Option  \text{II}:}$ $\x_{sol}=({\sum_{t=1}^{T}\gamma_t})^{-1}{\sum_{t=1}^{T}\gamma_t \hat\x_t}$ 
%%%%%%%%%%%%%%%%%%%%%%%
\end{algorithmic}%
\end{algorithm}
%%%%%%%%%%%%%%%%%%%%%%%ALG%%%%%%%%%%%%%%%%%%%%%%%
%%%%%%%%%%%%%%%%%%%%%%%%%%%%%%%%%%%%%%%%%%%%%%%%%

Next, we show that GG is able to find a $(1-1/e,\epsilon)$-approximate minimax solution after $\mathcal{O}(1/\epsilon^2)$ iterations. To prove this claim we require the following assumptions on the objective function $f$. 

%%%%%%%%%%%%%%%%%%%%%%%%
\begin{assumption}\label{assum:smooth}
The function $f$ is $L$-smooth with respect to $\x$, i.e., for any $\x,\x^{'}\in \mathbb{R}^d, S\in \mathcal{I}$, we have  $\|\nabla_{\x} f(\x,S)-\nabla_{\x} f(\x^{'},S)\|\leq L\|\x-\x^{'}\|$.
\end{assumption}

%%%%%%%%%%%%%%%%%%%%%%%%
\begin{assumption}\label{assum:bounded_grad}
The gradient of function $f$ with respect to $\x$ is uniformly bounded by a constant $M$, i.e.,  for any $\x\in \mathbb{R}^d, S\in 2^V$, we have  {$\|\nabla_\x f(\x,S)\|\leq M$}.
\end{assumption}

%%%%%%%%%%%%%%%%%%%%%%%%
\begin{theorem} \label{thm:GG2}
Consider Gradient Greedy (GG) in Algorithm~\ref{alg:(Extra)-gradient-greedy} {option \text{I}}. If $f$ is convex-submodular and Assumptions~\ref{assum:smooth}-\ref{assum:bounded_grad} hold, 
then the output of this algorithm after $\mathcal{O}(1/\epsilon^2)$ iterations  with step size $\mathcal{O}(\epsilon)$, is a $((1-{1}/{e}),\epsilon)$-approximate minimax solution of Problem~\eqref{eq:main_problem} . 
%%%%%%%%%%%%%%%%%%%%%%%
%\begin{align}\label{eq:GG main equation}
%({1-\frac{1}{e}})\max_S f({\x}_{sol},S) -\min_x\max_S f(x,S)\leq \frac{K}{\sqrt{T}}
% \end{align}
%
%%%%%%%%%%%%%%%%%%%%%%%
%%
\end{theorem}

%The result in Theorem~\ref{thm:GG2} shows that under smoothness (Assumption~\ref{assum:smooth}) and bounded gradient (Assumption~\ref{assum:bounded_grad}) Assumptions, Gradient-Greedy is able to achieve the optimal approximation guarantee of $(1-1/e)$ for the convex-submodular minimax problem in \eqref{eq:main_problem}.

The smoothness assumption (Assumption~\ref{assum:smooth}) is required to guarantee convergence of gradient-based methods at the rate of $1/\epsilon^2$. The bounded gradient assumption (Assumption~\ref{assum:bounded_grad}), however, comes from the fact that even in convex-concave problems gradient descent-ascent algorithms only converge when the gradient norm is uniformly bounded. This issue has been addressed in the convex-concave setting by the update of extra-gradient method which converges to a saddle point only under smoothness assumption. However, this improvement is not for free and it requires two gradient computations per update, instead of one. Next, we leverage this technique to present an alternating method that obtains the approximation factor and iteration complexity of GG without requiring Assumption~\ref{assum:bounded_grad}.  

%%%%%%%%%%%%%%%%%%%%%%%%%%%%%%%%%%%%%%%%%%
%\begin{wrapfigure}{R}{0.35\textwidth}
 %  \centering
  % \vspace{-6mm}
   % \includegraphics[width=1\linewidth]{Egradient} 
  %  \vspace{-8mm}
  %  \caption{Extra-gradient approach}
   % \vspace{-3mm}
%    \textcolor{red}{figure is good but too bid. Why do you have $X_t$? it should be $x_t$.}} \textcolor{blue}{because it is same notation as algorithms and other parts. I used $\x$. should I change all of them to $x$? }
  %  \label{fig:EGG} 
 % \end{wrapfigure}%%
  %%%%%%%%%%%%%%%%%%%%%%%%%%%%%%%%%%%%%%%%%%
   
    \vspace{2mm}
\noindent \textbf{Extra-gradient greedy.}
%\iffalse
%Before presenting our extra-gradient-based method, let us first briefly recap the update of extra-gradient for convex-concave minimax problems. If we consider $x_t$ and $y_t$ as the current minimization and maximization iterates, the main idea of Extra-gradient is to take a gradient step with respect to both minimization and maximization variables to find a set of auxiliary points $\hat{\x}_t$ and $\hat{\y}_t$, then find the new variables $\x_{t+1}$ and $\y_{t+1}$ by running  a gradient on the current iterates $\x_t$ and $\y_t$, where the gradients are evaluated at the auxiliary points $\hat{\x}_t$ and $\hat{\y}_t$. Specifically, the auxiliary points are computed as 
%\begin{align*}
 %   \hat{\x}_t= \pi_{X}(\x_{t} -\gamma_t  \nabla_x f(\x_{t},\y_t)), \qquad \text{and} \qquad   \hat{\y}_t= \pi_{Y}(\y_{t} -\gamma_t  \nabla_\y f(\x_{t},\y_t))
%\end{align*}
%%
%and they are used to compute the new variables $\x_{t+1}$ and $\y_{t+1}$ according to the updates:
%\begin{align*}
%   \x_{t+1} = \pi_{X}(\x_{t} -\gamma_t  \nabla_x f(\hat{\x}_t,\hat{\y}_t)), \qquad \text{and} \qquad  {\y}_{t+1}= \pi_{Y}(\y_{t} -\gamma_t  \nabla_\y f(\hat{\x}_t,\hat{\y}_t))
%\end{align*}
%%
% \fi
We now present the Extra-Gradient Greedy (EGG) algorithm, which consists of two gradient updates as suggested by extra-gradient and two greedy steps to find the auxiliary set $\hat{S_t}$ and the updated set $S_{t+1}$. In the extra-gradient method, we take a preliminary step to find a middle/auxiliary point and then compute the next iterate using the gradient information of the middle point. If we consider $\x_t$ and $S_t$ as the current iterates, we first run a gradient step to find the auxiliary minimization variable according to the update $\hat{\x}_t =  \pi_{\X}(\x_t -\gamma_t  \nabla f(\x_t,S_t))$
%\begin{equation*}
 %   \hat{\x}_t =  \pi_{\X}(\x_t -\gamma_t  \nabla f(\x_t,S_t))
%\end{equation*}
then we compute the auxiliary set $\hat{S_t}$ by performing a greedy step based on the auxiliary iterate $\hat{\x}_t$, i.e., $
\hat{S}_t = \textsc{Greedy}(f , k,\hat{\x}_t)
$. Once $\hat{\x}_t$ and $\hat{S_t}$ are computed,  we update the minimization variable $\x_{t+1}$ by descending towards a gradient evaluated at $(\hat{\x}_t,\hat{S}_t)$, i.e., $\x_{t+1} = \pi_{\X}( \x_t -\gamma_t  \nabla f(\hat{\x}_t,\hat{S}_t))$.
%\begin{equation*}
%\x_{t+1} = \pi_{\X}( \x_t -\gamma_t  \nabla f(\hat{\x}_t,\hat{S}_t))
%\end{equation*}
Lastly, we compute the new set $S_{t+1}$ by  running a greedy update based on the new iterate $\x_{t+1}$, i.e., $
S_{t+1} = \textsc{Greedy}(f , k,\x_{t+1})
$. Steps of EGG are outlined in Algorithm~\ref{alg:(Extra)-gradient-greedy} {(option \text{II})}.

Next we establish our theoretical result for Extra-gradient Greedy and show that only under smoothness assumption it finds an $(1-1/e,\epsilon)$-approximate minimax solution  after $\mathcal{O}(1/\epsilon^2)$ iterations.
%%%%%%%%%%%%%%%%%%%%%alg%%%%%%%%%%%%%%%%%%%%%%%%
\begin{algorithm}[t]
\caption{\\\text{$\textbf{Option  \text{I}:} $Gradient Replacement-greedy (GRG)}\\\text{$\textbf{Option  \text{II}:} $Extra-gradient Replacement-greedy (EGRG)}}\label{alg:(Extra)-gradient-Replacementgreedy}

\begin{algorithmic}
\State \textbf{Initialize} the set $S_1$ to $\emptyset$ and variable $\x_1$ to zero. 
\For{$t=1$ to $T$} 
%%%%%%%%%%%%%%%%%%%%%%%
 
%%%%%%%%%%%%%%%%%%%%%%%

   \State $\textbf{Option  \text{I}:}$ 
   $
{\x}_{t+1} =  \pi_{\X}(\x_{t} -\gamma_t  \nabla f(\x_{t},S_t))$
\State\qquad\qquad\,\,\,\,\,$
{S}_{t+1} = \textsc{RepGreedy}(f , k,{\x}_{t+1},S_{t})
$
\vspace{1mm}
%%%%%%%%%%%%%%%%%%%%%%%
\State$\textbf{Option  \text{II}:}$  
    $
\hat{\x}_t =  \pi_{\X}(\x_t -\gamma_t  \nabla f(\x_t,S_t))
$

\State\qquad\qquad\,\,\,\,\,\,\,\,$
\hat{S}_t = \textsc{RepGreedy}(f , k,{\x}_t,{S}_t)
$

\State\qquad\qquad\,\,\,\,\,\,\,\,$
\x_{t+1} = \pi_{\X}( \x_t -\gamma_t  \nabla f(\hat{\x}_t,\hat{S}_t))
$

\State\qquad\qquad\,\,\,\,\,\,\,\,$
{S}_{t+1} = \textsc{RepGreedy}(f , k,\hat{\x}_{t},\hat S_{t})
$
\vspace{-1mm}
\EndFor
%%%%%%%%%%%%%%%%%%%%%%%
\State $\textbf{Option  \text{I}:}$  ${\x}_{sol}=({\sum_{t=1}^{T}\gamma_t})^{-1}{\sum_{t=1}^{T}\gamma_t\x_t}$
\State$\textbf{Option  \text{II}:}$  $\x_{sol}=({\sum_{t=1}^{T}\gamma_t})^{-1}{\sum_{t=1}^{T}\gamma_t \hat\x_t}$ 
%%%%%%%%%%%%%%%%%%%%%%%
\end{algorithmic}%
\end{algorithm}
\vspace{-1mm}
%%%%%%%%%%%%%%%%%%%%%%%ALG%%%%%%%%%%%%%
\begin{theorem} \label{thm:EGG2}
Consider Extra-Gradient Greedy (EGG) outlined in Algorithm~\ref{alg:(Extra)-gradient-greedy} {option {II}}. If $f$ is convex-submodular and Assumption~\ref{assum:smooth} holds, then the output of this algorithm after $\mathcal{O}(1/\epsilon^2)$ iterations with step size $\mathcal{O}(\epsilon)$, is a $((1-{1}/{e}),\epsilon)$-approximate minimax solution of Problem~\eqref{eq:main_problem}. 
%%%%%%%%%%%%%%%%%%%%%%%
%\begin{align}\label{eq:EGG main equation}
%({1-\frac{1}{e}})\max_S f({\x}_{sol},S) -\min_x\max_S f(x,S)\leq \frac{K}{\sqrt{T}}
% \end{align}
%
%%%%%%%%%%%%%%%%%%%%%%%
\end{theorem}

\begin{remark}
Note that as both GG and EGG are greedy based methods, they can also be used for the case of general matroid constraint. However, the approximation guarantee would be $1/2$ instead of $1-1/e$. The details are provided in the Appendix.
\end{remark}

 \subsection{Replacement Greedy-based Methods}
As we showed earlier, for the cardinality constraint problem GG and EGG achieve the optimal approximation guarantee of $1-1/e$ for the minimax problem in \eqref{eq:main_problem}. However, they both require running greedy updates at each iteration which makes their per iteration complexity $\mathcal{O}(nk)$. To resolve this issue, we propose the use of replacement-greedy in lieu of greedy update. This modification reduces the complexity of each iteration to $\mathcal{O}(n+k)$ at the cost of lowering the approximation factor.
 
 \vspace{2mm}
\noindent{\textbf{Gradient replacement-greedy.}}
We first present the Gradient Replacement-Greedy (GRG) algorithm which alternates between a gradient update and a replacement greedy update. As shown in  Algorithm~\ref{alg:(Extra)-gradient-Replacementgreedy} {option \text{I}}, the only difference between GRG and GG algorithms is the substitution of greedy update with replacement greedy. Next, we establish the theoretical guarantee of GRG.
 %\iffalse
    %%%%%%%%%%%%%%%%%%%%%alg%%%%%%%%%%%%%%%%%%%%%%%%
%begin{algorithm}[t]
%\caption{\text{Gradient Replacement-greedy(GRG)}}\label{alg:Gradient Replacement-greedy}
%\begin{algorithmic}
% \State \textbf{Initialize} the set $S_1$ to the empty set and variable $x_1$ to the zero. 
%\For{$t=1$ to $T$} 
 %\vspace{2mm}
%%%%%%%%%%%%%%%%%%%%%%%
%     \State   Update $\x_{t+1}$ using gradient descent  $
%{\x}_{t+1} =  \pi_{X}(\x_{t} -\gamma_t  \nabla f(\x_{t},S_{t}))
%$
% \vspace{2mm}
%%%%%%%%%%%%%%%%%%%%%%%
%    \State Update $ S_{t+1}$ using replacement greedy update ${S}_{t+1} = \textsc{RepGreedy}(f %, k,{\x}_{t+1},S_{t})
%$
   %%%%%%%%%%%%%%%%%%%%%%%
 %   % \STATE  \textcolor{Brown}{$$y^{*}=\arg\max\limits_{e\in V}\{f( \hat \x_t,\hat %S_{t-1}\setminus \{y\})\}$$}
    %\vspace*{-5mm}
%%%%%%%%%%%%%%%%%%%%%%%
%\EndFor
%%%%%%%%%%%%%%%%%%%%%%%
%\State Return solution $\x_{sol}=({\sum_{t=1}^{T}\gamma_t})^{-1}{\sum_{t=1}^{T}\gamma_t \x_t}$
%\end{algorithmic}%
%\end{algorithm}
%%%%%%%%%%%%%%%%%%%%%%%%%%%%%%%%%%%%%%%%%%%%%
%\fi

\begin{theorem} \label{thm:GRG2}
Consider the {Gradient Replacement-Greedy (GRG)} algorithm in Algorithm~\ref{alg:(Extra)-gradient-Replacementgreedy} {option \text{I}}. If $f$ is convex-submodular and Assumptions~\ref{assum:smooth}-\ref{assum:bounded_grad} hold, then the output of this algorithm after $\mathcal{O}(1/\epsilon^2)$ iterations with step size $\mathcal{O}(\epsilon)$, is a $(1/2,\epsilon)$-approximate minimax solution of Problem~\eqref{eq:main_problem}. 
%Consider the $\textbf{Gradient Replacement-greedy(GRG)}$ algorithm outlined in Algorithm~\ref{alg:Gradient Replacement-greedy}, if the functions $f$ is convex-submodular and $G-$gradient bounded , then the point ${\x}_{sol}=\frac{\sum_{t=1}^{T}\gamma_t\x_t}{\sum_{t=1}^{T}\gamma_t}$  gives us $(\frac{1}{2+\frac{k}{k-1}},\frac{K}{\sqrt{T}})$-solution for minimax convex-submodular problem on cardinality constraint where $\gamma_t=\frac{c}{\sqrt{T}}$ and $c,K$ are some constants:
%%%%%%%%%%%%%%%%%%%%%%%%
%\begin{align}\label{eq:EGG main equation}
%\frac{1}{2+\frac{k}{k-1}}\max_S f({\x}_{sol},S) -\min_x\max_S f(x,S)\leq \frac{K}{\sqrt{T}}
% \end{align}
%%%%%%%%%%%%%%%%%%%%%%%
\end{theorem}
%The proof of theorem~\ref{thm:GRG2} can be find in supplementary material.

\noindent{\textbf{Extra-gradient replacement-greedy.}} The GRG algorithm requires the bounded gradient assumption similar to GG. To address this issue, a natural idea is to exploiting the  extra-gradient approach for updating $\x$ and introducing the Extra-gradient Replacement-greedy (EGRG) algorithm, outlined in Algorithm~\ref{alg:(Extra)-gradient-Replacementgreedy} {option \text{II}}.
However, unlike the case of Greedy-based methods, here we can not drop the bounded gradient assumption by exploiting the idea of extra-gradient update. Next, we elaborate on this issue.

Note that to prove that EGRG finds a $({1}/{2},\epsilon)$-approximate minimax solution we need to find an upper bound on 
$f(\hat{\x}_t,S)-2f(\hat{\x}_t,\hat {S}_t)$ for every $S$. To establish such a bound, we need to relate $f(\hat{\x}_t,\hat {S}_t)$ to $f({\x}_t,\hat{S}_t)$ which requires the function $f$ to be Lipschitz with respect to $\x$, which is equivalent to the bounded gradient condition in  Assumption~\ref{assum:bounded_grad}; see proof of Theorem 2 in the appendix for more details. Note that such argument is not required for the EGG method, as in greedy based method we always have the following  inequality $f(\hat{\x}_{t},S)-(1-{1}/{e})^{-1}f(\hat{\x}_{t},\hat{S}_{t})\leq0$ for every $S$. As a result, the required conditions for the convergence of GRG and EGRG are similar and we only state EGRG results for completeness. 

%Each step consists of Extra-gradient step and replacement-greedy step. In the Extra-gradient step, similar to EGG, we update $\x_t$ to $\hat \x_{t}$ using projected gradient decent $\hat{\x}_t =  \pi_{X}(\x_t -\gamma_t  \nabla f(\x_t,S_t))$ and then it will use the replacement-greedy to find  $\hat S_t$. Replacement-greedy replaces the worst element in the current set with best element in the ground set. Replacement-greedy will find an  element $\hat e^{*}$ that replacing it with its replacement $\text{REP}(\hat e^{*},\hat S_{t-1},\hat \x_t)$ at point $\hat \x_t$ results in maximum gain. In second step, it will update $\x_{t+1}$ with projected gradient descent using gradient at point $(\hat{\x}_t,\hat{S}_t)$ as  $\x_{t+1} = \pi_{X}( \x_t -\gamma_t  \nabla f(\hat{\x}_t,\hat{S}_t))$. Then it computes the element $ e^{*}$ in the ground set which has maximum replacement gain as  $e^{*}= \arg\max\limits_{e\in V}\{\nabla_e f(\x_{t+1},\hat S_{t})\}$ to construct $S_{t+1}$ as $S_{t+1}= \hat{S}_{t}\cup \{\hat e^*\}\setminus \{\text{REP}( e^{*},\hat S_{t},\x_{t+1})\}$. We illustrate EGRG algorithm in the figure \ref{fig:EGRG}.

\begin{theorem} \label{thm:EGRG2}

Consider {Extra-Gradient Replacement Greedy(EGRG)} in Algorithm~\ref{alg:(Extra)-gradient-Replacementgreedy} {option \text{II}}. If $f$ is convex-submodular and Assumptions~\ref{assum:smooth}-\ref{assum:bounded_grad} hold, then the output of EGRG after $\mathcal{O}(1/\epsilon^2)$ iterations with step size $\mathcal{O}(\epsilon)$, is a $({1}/{2},\epsilon)$-approximate minimax solution of \eqref{eq:main_problem}. 
%Consider the $\textbf{Extra-Gradient Replacement Greedy(EGRG)}$ algorithm outlined in Algorithm~\ref{alg:Extra-gradient-replacement-greedy}, if the functions $f$ is convex-submodular, $G-$gradient bounded, and L-smooth then the point ${\x}_{sol}=\frac{\sum_{t=1}^{T}\gamma_t\hat\x_t}{\sum_{t=1}^{T}\gamma_t}$  gives us $(\frac{1}{2},\frac{K}{\sqrt{T}})$-solution for minimax convex-submodular problem on cardinality constraint where $\gamma_t=\frac{c}{\sqrt{T}}$ and $c,K$ are some constants:
%%%%%%%%%%%%%%%%%%%%%%%
%\begin{align}\label{eq:EGG main equation}
%\frac{1}{2}\max_S f({\x}_{sol},S) -\min_x\max_S f(x,S)\leq \frac{K}{\sqrt{T}}
% \end{align}
%%%%%%%%%%%%%%%%%%%%%%%
\end{theorem}
%
%The proof of theorem~\ref{thm:EGRG2} can be find in supplementary material.

%It is worth noting this improvement comes at the cost of computing two gradients and running two replacement greedy updates per iteration.  
\vspace{-1mm}

 \subsection{Extra-gradient on Continuous Extension}

  So far all proposed algorithms achieve $(\alpha,\epsilon)$-approximate minimax solutions in $\mathcal{O}({1}/{\epsilon^2})$ iterations. In this section, we investigate the possibility of achieving a faster rate of  $ \mathcal{O}({1}/{\epsilon})$. Note that, in the discussed algorithms, the update for the discrete variable is not smooth and the iterates jump from one set to another in consecutive iterations, which results in slowing down the convergence. To overcome this limitation, we introduce the continuous multi-linear extension of Problem~\eqref{eq:main_problem}; for introduction to multi-linear extension of submodular maximization problems and how to optimize them see \citep{calinescu2011maximizing,badanidiyuru2014fast,feldman2011unified, hassani2017gradient,mokhtari2020stochastic,hassani2020stochastic, sadeghi2020online}. As we will show, the continuous extension of Problem~\eqref{eq:main_problem} is equivalent to its original version, and by extending the extra-gradient methodology to this setting we achieve a convergence rate of $\mathcal{O}({1}/{\epsilon})$ for the case that $\mathcal{I}$ is a~matroid. \looseness=-1
%{\defi Let $\I$ be a nonempty family of allowable subsets of the ground set $V$, then the tuple $(V, \I)$ is a \textit{matroid} if and only if the following conditions hold (i)  For any $A\subset B\subset V$, if $B\in\I$, then $A\in\I$ (ii) For all $A, B\in\I$, if $|A|<|B|$, then there is an $e\in B\backslash A$ such that $A \cup\{e\} \in\I$.
%}

  %%%%%%%%%%%%%%%%%%%%%alg%%%%%%%%%%%%%%%%%%%%%%%%
\begin{algorithm}[t]
\caption{\text{Extra-gradient on Continuous Extension}}\label{alg:Extra-gradient on Continuous Extension}
\begin{algorithmic}
\State \textbf{Initialize} the variables $\y_1$ and $\x_1$ to zero. 
\For{$t=1$ to $T$} 
% \vspace{1mm}
%%%%%%%%%%%%%%%%%%%%%%%
    \State$\hat{\x}_t =  \pi_{\X}(\x_t -\gamma_t  \nabla_x F(\x_t,\y_t))$
 \vspace{1mm}
%%%%%%%%%%%%%%%%%%%%%%%
   \State$
\hat{\y}_t =  \pi_{\K}(\y_t +\gamma_t  \nabla_y F(\x_t,\y_t))
$
%%%%%%%%%%%%%%%%%%%%%%%
   %\STATE  \textcolor{Brown}{$$e^{*}=\arg\max\limits_{e\in V}\{\Delta_e f( \hat \x_t,\hat S_{t-1}\setminus \{y^{*}\})\}$$}
   % \vspace*{-5mm}
%%%%%%%%%%%%%%%%%%%%%%%
    %  \STATE  \textcolor{Brown}{ $$\hat{S}_{t}=\hat S_{t-1}\cup \{e^*\}\setminus \{y^{*}\} $$}
    %   \vspace*{-5mm}
%%%%%%%%%%%%%%%%%%%%%%%
    \vspace{1mm}
%%%%%%%%%%%%%%%%%%%%%%%
    \State$
\x_{t+1} = \pi_{\X}( \x_t -\gamma_t  \nabla_x F(\hat{\x}_t,\hat{\y}_t))
$
 \vspace{1mm}
%%%%%%%%%%%%%%%%%%%%%%%
\State$ 
\y_{t+1} = \pi_{\K}( \y_t +\gamma_t  \nabla_y F(\hat{\x}_t,\hat{\y}_t))
$
%%%%%%%%%%%%%%%%%%%%%%%
\EndFor
%%%%%%%%%%%%%%%%%%%%%%%
\State Return solution $\x_{sol}=({\sum_{t=1}^{T}\gamma_t})^{-1}{\sum_{t=1}^{T}\gamma_t\hat \x_t}$
\end{algorithmic}%
\end{algorithm}

  {\defi The \emph{continuous extension} of a function $f: \mathbb{R}^d \times 2^V \to \mathbb{R}_{+}$ is the function $F: \mathbb{R}^d \times [0,1]^n \to \mathbb{R}_{+}$ defined as 
    $F(\x,\y) = \mathbb{E}_{S \sim \y} [f(\x,S)],$ 
where  $S\sim \y$ is a random set wherein each element $i$ is included with probability $y_i$ independently.}

We show that for convex-submodular problems we have (see Proposition~\ref{prop:1} in Appendix~\ref{apendix:Extra Gradient on continuous Extension}):  
\begin{equation}
    \min_{\x \in \mathcal{X}} \max_{S \in \mathcal{I}} f(\x,S) = \min_{\x\in \mathcal{X} } \max_{\y \in \mathcal{K}} F(\x,\y),
\end{equation}
where $\mathcal{I}$ is assumed to be a matroid constraint and $\mathcal{K}$ is the corresponding base polytope($\mathcal{K}=\textbf{conv}\{1_{S} : S\in\mathcal{I} \}$).  
  We present Extra-Gradient on Continuous Extension  (EGCE) in Algorithm~\ref{alg:Extra-gradient on Continuous Extension} which applies the updates of extra-gradient on the continuous extension  function $F(\x,\y)$.

  %This iterative algorithm  first perform projected gradient descent on $\x$ and 
  %$\y$ to obtain middle point $(\hat \x_t, \hat \y_t)$ as $$\hat{\x}_t =  \pi_{X}(\x_t -\gamma_t  \nabla_x F(\x_t,\y_t)), \qquad 
%\hat{\y}_t =  \pi_{X}(\y_t -\gamma_t  \nabla_y F(\x_t,\y_t))
%$$
%Then, it uses gradient at the middle point to update $\x$ and $\y$ as 
%$$
%\y_{t+1} = \pi_{X}( \y_t -\gamma_t  \nabla_y %F(\hat{\x}_t,\hat{\y}_t)),\qquad
%\x_{t+1} = \pi_{X}( \x_t -\gamma_t  \nabla_x F(\hat{\x}_t,\hat{\y}_t)).
%$$

 \begin{theorem} \label{thm:EGcont2}
 Consider the {Extra-Gradient On Continuous Extension} (EGCE) algorithm outlined in Algorithm~\ref{alg:Extra-gradient on Continuous Extension}. If $f$ is convex-submodular and Assumptions~\ref{assum:smooth}-\ref{assum:bounded_grad} hold, then the output of this algorithm  after $\mathcal{O}(1/\epsilon)$ iterations with constant step size is a $(1/2,\epsilon)$-approximate minimax solution of Problem~\eqref{eq:main_problem}.
 %
%Consider the $\textbf{Extra-Gradient On Continuous Extension}$ algorithm outlined in Algorithm~\ref{alg:Extra-gradient on Continuous Extension}, if the functions $f$ is convex-submodular , then the point ${\x}_{sol}=\frac{\sum_{t=1}^{T}\gamma_t\hat\x_t}{\sum_{t=1}^{T}\gamma_t}$  gives us $(0.5,\frac{K}{{T}})$-solution for minimax convex-submodular problem on matroid constraint where $\gamma_t=\frac{c}{\sqrt{T}}$ and $c,K$ are some constants:
%%%%%%%%%%%%%%%%%%%%%%%
%\begin{align}\label{eq:EGcont main equation}
%\frac{1}{2}\max_S f({\x}_{sol},S) -\min_x\max_S f(x,S)\leq \frac{K}{{T}}
% \end{align}
 %
%%%%%%%%%%%%%%%%%%%%%%%
%
\end{theorem}

\section{Experiments}\label{sec:Poisoning Attack}
  In this section, we study two specific instances of Problem~\eqref{eq:main_problem}: (i) convex-facility location functions along with a synthetic experimental setup, and (ii) designing adversarial attacks for item recommendation which is a real world application of our framework.   

\textbf{Convex-facility location functions.}
Consider the function $f: \mathbb{R}^d \times 2^V \to \mathbb{R}_{+}$ defined as
%\begin{equation}\label{eq:Convex-Facility Location Function}
    $f(\x,S) = \sum_{i = 1}^n \max_{j \in S} f_{i,j}(\x) + g(\x)$,
%\end{equation}
where $g: \mathbb{R}^d \to \mathbb{R}$ and $f_{i,j}: \mathbb{R}^d \to \mathbb{R} $ are convex. Indeed, $f(\x,S)$ is convex with respect to $\x$. Also, for a fixed $\x$, we recover the objective of the facility location problem, which is submodular and monotone.  
%, and the objective function in \eqref{eq:Convex-Facility Location Function} is a special case of Problem~\eqref{eq:main_problem}.
To introduce our setup, suppose  $\x\in \mathbb{R}_+^d$ can be written as the concatenation of $n=d/m$ vectors $\x_i\in \mathbb{R}_{+}^m $ of size $m$, i.e., $\x=[\x_1;\dots;\x_n]$.
In our experiments, we assume that the function $f_{i,j}(\x)$ is defined as $f_{i,j}(\x)=\x_i^TQ_{i,j}\x_j$, where $Q_{i,j}\in \mathcal{S}_{++}^{m}$ is a positive definite matrix and all of its elements are also positive, i.e.,  $Q_{i,j}>0$. Moreover, we consider the case that the regularization  function $g$ is defined as  $g(\x):={\lambda}(\sum_{i=1}^n\|\x_i\|^2)^{-1}$, and the constraint set for the minimization variable $\x$ is defined as $\mathcal{C}:=\{\x=[\x_1;\dots;\x_n]| \|\x_i\|\leq 1,\ \text{for}\ i=1,\dots,n \}$. Considering these definitions the convex-submodular minimax optimization problem that we aim to solve can be written as
 %Let us define for $i\in[n]$, $\x_i\in \mathcal{C}$  where set $\mathcal{C}$ is $\mathcal{C}:=\{\x|\x\in \mathbb{R}_{+}^{m},\|\x\|\leq 1\}$. Next, we define $\x$ as collection of all the vectors $\x_i$ for $i\in[n]$, $\x=vec(\x_1,\x_2,\dots, \x_n)$. Then we let $f_{i,j}(\x)=\x_i^TQ_{i,j}\x_j$ positive definite, entry-wise positive matrix $Q_{i,j}$ (i.e. $Q_{i,j}\in \mathcal{S}_{++}^{m}, Q_{i,j}>0$), and function $g(\x):={\lambda}(\sum_{i=1}^n\|\x_i\|^2)^{-1}$.
% We examine our algorithms for problem in \eqref{eq:simulation}
\begin{equation}\label{eq:simulation}
\min_{\x_i\in \mathcal{C}}\max_{|S|\leq k} \sum_{i = 1}^n \max_{j \in S} \x_i^TQ_{i,j}\x_j+\lambda
   \Big(\sum\limits_{i=1}^n\|\x_i\|^2\Big)^{-1}
\end{equation}
where the constraint on the maximization variable $S$ is a cardinality constraint of size $k$.
For our numerical experiments, we tested two cases, in the first case we set the problem parameters as $m=10$, $n=30$, $k=5$, and $\lambda=1$ and  in the second case we set the problem parameters as $m=10$, $n=50$, $k=10$, and $\lambda=1$.

  \begin{figure}[t]

\centering
\begin{minipage}{0.40\textwidth}
     \begin{figure}[H]
        \includegraphics[width=\linewidth]{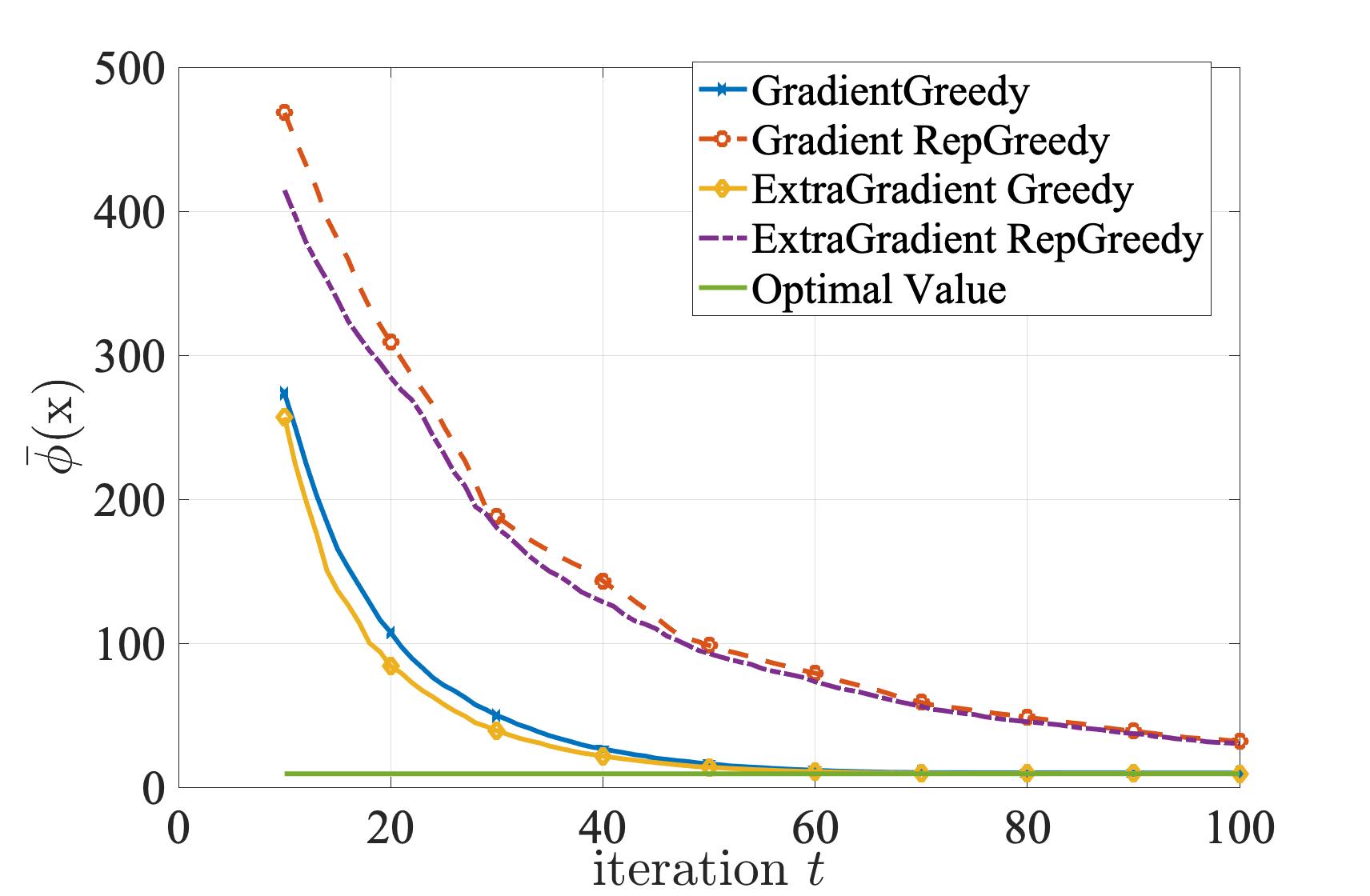}
        \label{fig:syn1}  
        %\caption{Comparison of Gradient Greedy,Gradient Replacement Greedy(Gradient RepGreedy), Extra gradient Greedy, and Replacement Greedy(Extra gradient RepGreedy) for Synthetic example in \eqref{eq:simulation}.}
    \end{figure}
    \end{minipage}
    \hspace{-4mm}
\begin{minipage}{0.40\textwidth}
    \begin{figure}[H]
        \includegraphics[width=\linewidth]{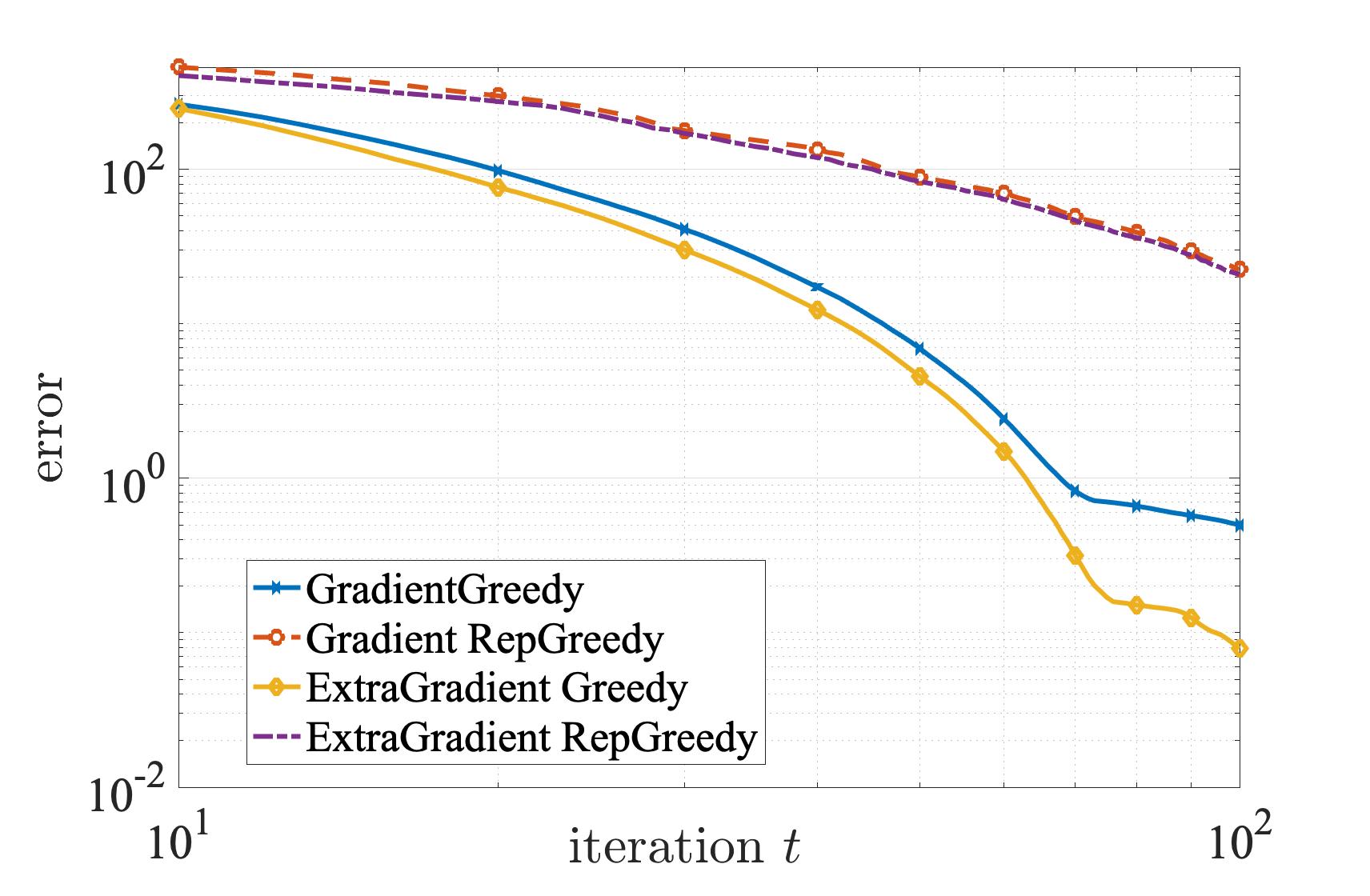}
        \label{fig:synlog2}
        %\caption{Comparison of Gradient Greedy,Gradient Replacement Greedy(Gradient RepGreedy), Extra gradient Greedy, and Replacement Greedy(Extra gradient RepGreedy) for Synthetic example in \eqref{eq:simulation}.}
    \end{figure} %
    \end{minipage}
    \hspace{-4mm}
    \vspace{-3mm}
    
    \vspace{-5mm}

\begin{minipage}{0.40\textwidth}
     \begin{figure}[H]
        \includegraphics[width=\linewidth]{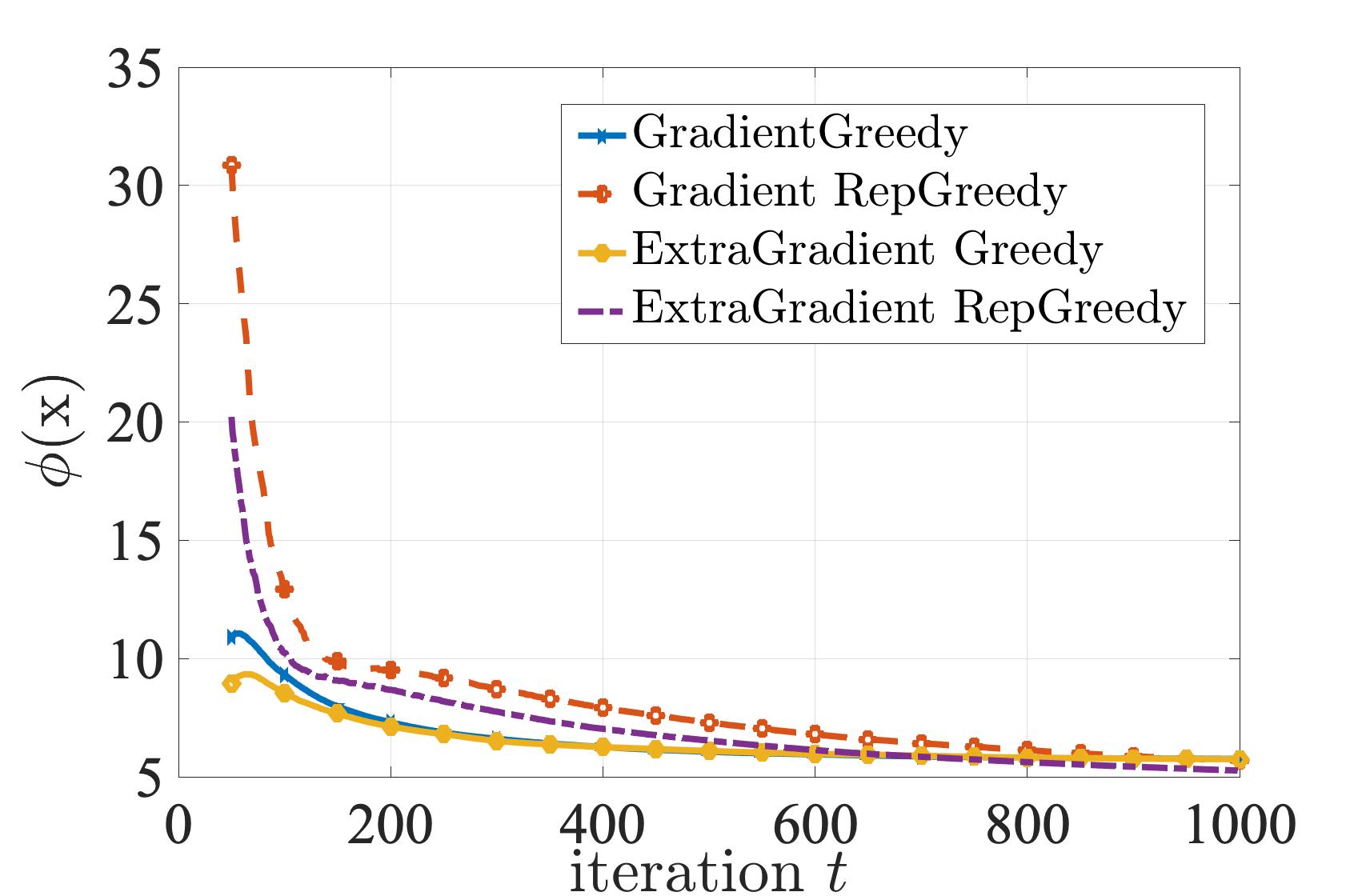}
        \label{fig:syn}  
        %\caption{Comparison of Gradient Greedy,Gradient Replacement Greedy(Gradient RepGreedy), Extra gradient Greedy, and Replacement Greedy(Extra gradient RepGreedy) for Synthetic example in \eqref{eq:simulation}.}
    \end{figure}
    \end{minipage}
    \hspace{-2mm}
\begin{minipage}{0.40\textwidth}
    \begin{figure}[H]
        \includegraphics[width=\linewidth]{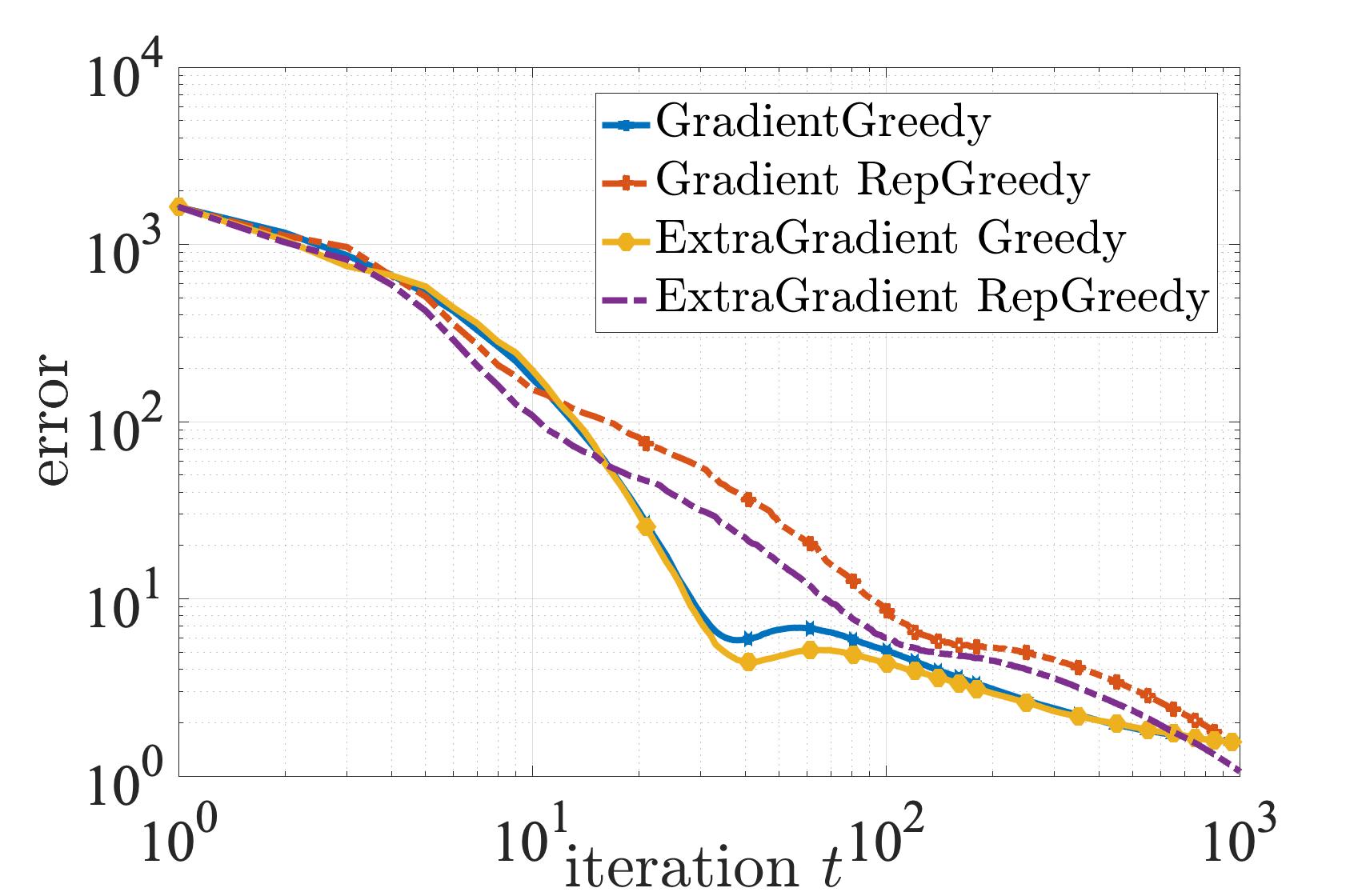}
        \label{fig:synlog}
        %\caption{Comparison of Gradient Greedy,Gradient Replacement Greedy(Gradient RepGreedy), Extra gradient Greedy, and Replacement Greedy(Extra gradient RepGreedy) for Synthetic example in \eqref{eq:simulation}.}
    \end{figure} %
    \end{minipage}
    \hspace{-2mm}

    \vspace{-3mm}
    
    \caption{Comparison of our proposed methods for  convex-facility location functions(case I  and Case II)}\label{fig:confac_opt}
    \vspace{-2mm}
\end{figure}

\underline{Case I ($m=10$, $n=30$, $k=5$, $\lambda=1$).}
In this case, we choose $m,n$ to be small so that we can solve the inner max in problem \eqref{eq:simulation} and compute $\bar\phi(\x)=\max_{|S|\leq k}f(\x,S)$ exactly using search over all the subsets of size $k$. We report $\bar\phi(\x_t)$ as well as optimal value of problem \eqref{eq:simulation}. Results in Figure \ref{fig:confac_opt}. (first plot) show that the algorithms converge to the optimal minimax value. We also demonstrate the relative error of these algorithms $\text{error}_t:=\phi(\x_t)-\text{OPT}$ in second plot. As we observe in Figure \ref{fig:confac_opt} (second plot), greedy based methods converge faster than replacement greedy based algorithms in terms of iteration complexity. %The second plot in Figure  \ref{fig:confac_opt} is in log-log scale. %In this plot, the greedy based methods have slope of -1, and the replacement greedy based methods have slope of -0.5 which indicates convergence rates better than $\mathcal{O}(1/\sqrt{t})$. These convergence rates are better than worst-case bounds in our theoretical results.

\underline{Case II ($m=10$, $n=50$, $k=10$, $\lambda=1$).}
We now investigate the behavior of our proposed methods for solving  \eqref{eq:simulation} in the second case when $k,n$ are relatively larger. Note that exact computation of  $\bar\phi(\x)=\max_{|S|\leq k}f(\x,S)$ is not computationally tractable for this case, since it requires solving a submodular maximization problem. Hence, in third plot in figure \ref{fig:confac_opt}, we report the value of the function  $\phi(\x):=f(\x,\textsc{Greedy}(f,k,\x))$ which is an approximation  for $\bar\phi(\x)$. In other words, instead of computing $\bar\phi(\x)$ which is the maximum of $f(\x,S)$ over the choice of $S$, we report $\phi(\x)$ which is the value of $f(\x,S)$ when $S$ is obtained via the greedy method.
The convergence paths of $\phi(\x_t)$ for our proposed methods are reported in the third plot of Figure~\ref{fig:confac_opt}. We further show the relative error of these algorithms defined as $\text{error}_t:=\phi(\x_t)-\phi(\x_{T})$ in the fourth plot to better compare their convergence rates.
%As we observe in the second plot, the convergence slope of all considered methods is almost $-1$ in the log-log scale, which indicates a convergence rate of $\mathcal{O}(1/t)$. Hence, to achieve an $\epsilon$-accurate solution we require $1/\epsilon$ iterations, which is better than our worst-case theoretical results.
%%%%%%%%%%%%%%%%%%%%%%%%%%%%%

%%%%%%%%%%
  \begin{figure}[t]
   \vspace{-1mm}
     \begin{minipage}{0.5\textwidth}
     \begin{figure}[H]
        \includegraphics[width=\linewidth]{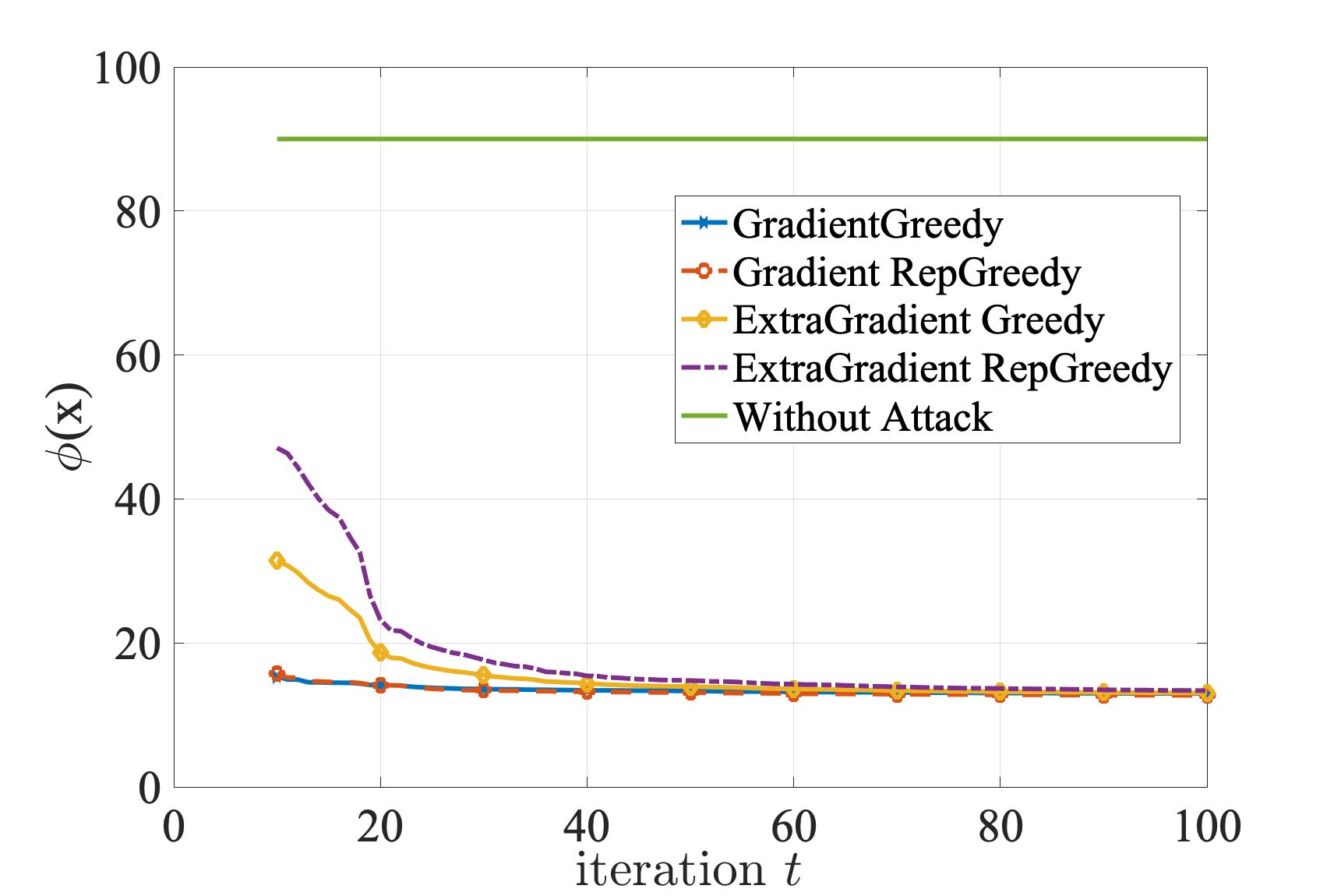}
        \label{fig:syn3}  
        %\caption{Comparison of Gradient Greedy,Gradient Replacement Greedy(Gradient RepGreedy), Extra gradient Greedy, and Replacement Greedy(Extra gradient RepGreedy) for Synthetic example in \eqref{eq:simulation}.}
    \end{figure}
    \end{minipage}
     \hspace{-4mm}
     \vspace{-1mm}
    \begin{minipage}{0.5\textwidth}
     \begin{figure}[H]
        \includegraphics[width=\linewidth]{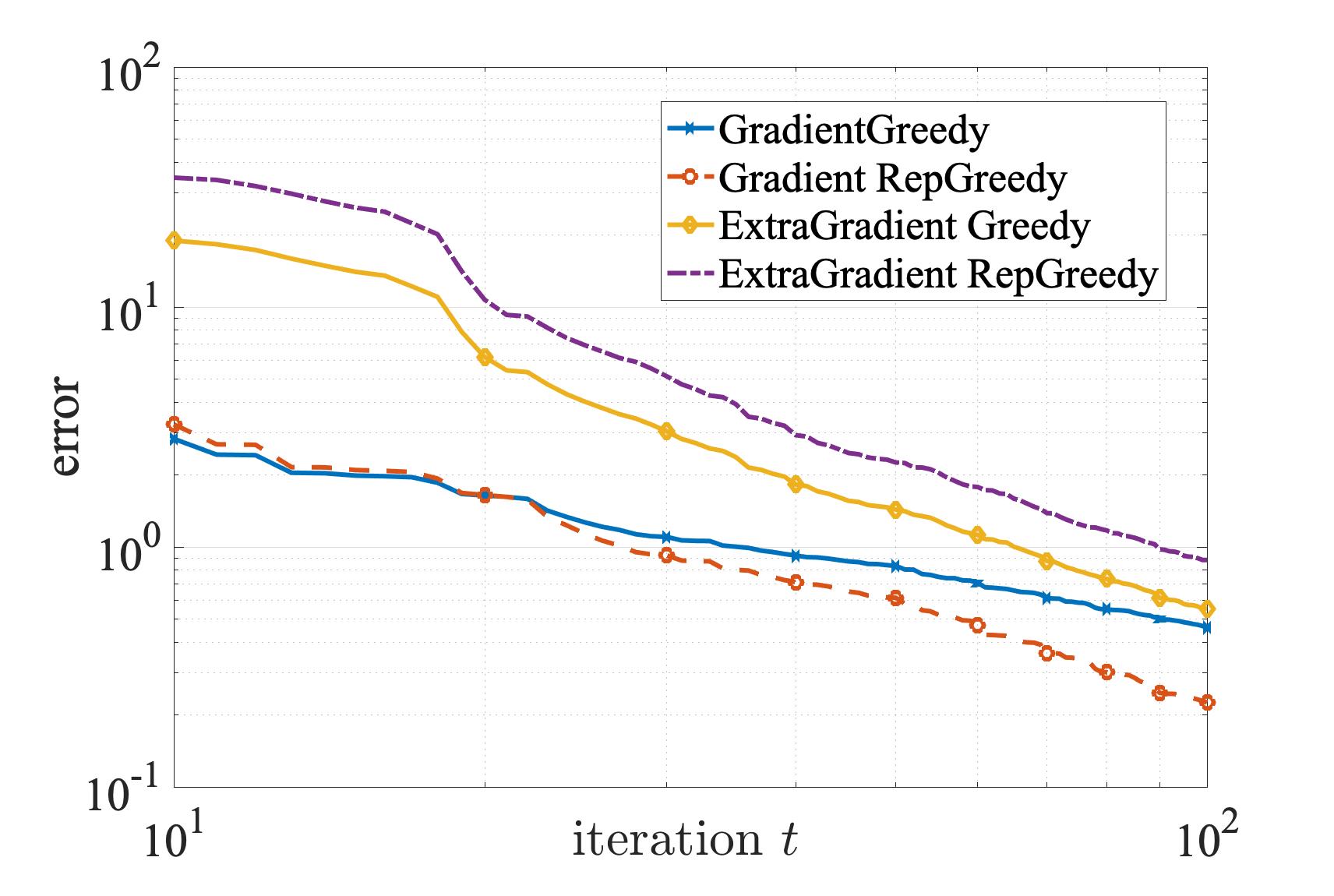}
        \label{fig:syn4}  
        %\caption{Comparison of Gradient Greedy,Gradient Replacement Greedy(Gradient RepGreedy), Extra gradient Greedy, and Replacement Greedy(Extra gradient RepGreedy) for Synthetic example in \eqref{eq:simulation}.}
    \end{figure}
    \end{minipage}
    \vspace{-4mm}
    \vspace{-1mm}
    \caption{Comparison of our proposed methods for  for Problem \eqref{eq:mainadvmov} (the green line is the performance of the recommender system when there is no adversary)}\label{fig:advatt}
    \vspace{-2mm}
\end{figure}

\textbf{Adversarial Attack for Item Recommendation.} In this section, we study the application of designing an adversarial attack for a movie recommendation task.
  Consider a (completed) rating matrix $X$ whose entries $X_{i,j}$ correspond to the estimated rating that user $i$ has given to  movie $j$. Given a  rating matrix $X$, the recommender system chooses $k$ movies  via maximizing the utility function  $\max_{|S| \leq k} h(X,S):=\frac{1}{|\mathcal{U}|}\sum_{u\in\mathcal{U}}\max_{j\in S} X_{u,j} $
where $\mathcal{U}$ is the set of all users.
The attacker's goal is to slightly  perturb the rating matrix $X$ to a matrix $X'$ such that the utility $\max_{|S| \leq l} h(X', S)$ is minimized. 
Therefore, the attacker aims at solving the minimax problem 
\begin{equation}\label{eq:mainadvmov}
\min\limits_{\substack{\|X'-X\|_{F}\leq \epsilon\\0\leq X'_{i,j}\leq 5} } \max\limits_{|S| \leq k} h(X',S) ,
\end{equation}
 where $\|.\|_{F}$ is the Frobenius norm. Note that $h(X,S)$ is convex-submodular (convexity in $\x$ is clear, and the function $h(\x, S )$ is a facility location function in $S$). Hence, this problem is an instance of Problem~\eqref{eq:main_problem}.
To evaluate the performance of our methods, we consider movie recommendation on the Movielens dataset \citep{harper2015movielens}.  We pick 2000 most rated movies with 200 users with highest number of rates for these movies (similar to \citep{stan2017probabilistic,adibi2020submodular}) and we set $k=10$. The adversary has a power to manipulate up to $0.5\%$ of  movies ratings on average (i.e. $\epsilon=0.5\times 0.01\times 200\times 2000$). We plot $\phi(X_\text{alg})$ in each iteration as a measure of effectiveness of our algorithms and compare it to the case that there is no attack. Figure \ref{fig:advatt} shows the comparison of our algorithms.
As we can see in Figure \ref{fig:advatt}, the facility location based recommendation systems are extremely vulnerable to adversarial attacks and the performance drops from 90 (when there is no adversary) to around 12 when we have attacks.  
 %First plot of Figure \ref{fig:advatt} depicts convergence path of $\phi(X_{t})$ and the second plot
 %reports $\text{error}_t=\phi(X_t)-\phi(X_{T})$. Both replacement greedy based algorithms and greedy based algorithms have a slope close to -1/2 in log-log scale which indicates a convergence rate of $\mathcal{O}(1/\sqrt{t})$, meaning we need $\mathcal{O}(1/\epsilon^2)$ iterations to achieve $\epsilon$ accuracy, which matches our theoretical bounds in Section~\ref{sec:algs}.

\section{Conclusion}

In this paper, we introduced and studied the convex-submodular minimax problem in \eqref{eq:main_problem}. 
We defined multiple notions of (near-) optimality and provided hardness results regarding these notions in various regimes. In particular, one of the notions was  $(\alpha,\epsilon)$-approximate minimax solution. We showed that for $\alpha>1-1/e$ finding an $(\alpha,\epsilon)$-approximate minimax solution is hard.  
For $\alpha\leq 1-1/e$, we proposed five algorithms and characterized their theoretical guarantees in different settings. The main take-away message from our algorithmic procedures is that, if the function $f$ has bounded gradient, then one can use the GG Algorithm, or the GRG algorithm which has a better complexity albeit it has a worse approximation factor. If the gradient of $f$ is not uniformly bounded, then one has to resort to the proposed EGG algorithm. %We refer to Table~\ref{tab:1} for the complexity and solution quality of these methods. 

\section*{Acknowledgement}
The work of A. Adibi and H. Hassani is funded by NSF award CPS-1837253, NSF
CAREER award CIF-1943064, and Air Force Office of Scientific Research Young Investigator Program (AFOSR-YIP) under award FA9550-20-1-0111. The work of A. Mokhtari is supported in part by NSF Grant 2007668, ARO Grant W911NF2110226, the Machine Learning Laboratory at UT Austin, and the NSF AI Institute for Foundations of Machine Learning.

%%%%%%%%%%%%%%%%%%%%%%%%%%%%%
%\bibliographystyle{ieeetr}

\bibliography{example_paper,ref}
%%%%%%%%%%%%%%%%%%%%%%%%%%%%%

%%%%%%%%%%%%%%%%%%%%%%%%%%%%%

 \appendix
\onecolumn

\section{Appendix}
%%%%%%%%%%%%%%%%%%%%%%%%%%%%%

\subsection{Proof of Theorem \ref{thm:saddle_neg}}\label{sec:proof of negative result1}
%%%%%%%%%%%%%%%%%%%%%%%%%%%%
%%%%%%%%%%%%%%%%%%%%%%%%%%%%
Consider the function $f:\mathbb{R}^d \times 2^V \to \mathbb{R}_{+}$, where $f(\x,.)$ is submodular for every $\x$ and $f(.,S)$ is convex for every $S$. Then, the maxmin convex-submodular problem is an optimization problem where the maximization is over continuous variable and minimization is over a discrete variable as
\begin{equation} \label{eq:main_problem22}
   {\rm{OPT}}_{maxmin} \triangleq \max_{S \in \mathcal{I}}\min_{x \in \mathcal{X}} f(\x,S),
\end{equation}
 
Let us define the notion of approximate solution for maxmin problem as follows:
{\defi We call a point $\hat{S}$ an $({\alpha},\epsilon)$-approximate maxmin solution of Problem~\eqref{eq:main_problem22} if it satisfies
%%%%%%%%%%%%%%%%%%%%%%%%%%%%
\begin{equation}\label{eq:(alpha,eps)solx}
 {\alpha}\underline{\phi}(\hat{S})\geq  {\rm{OPT}_{maxmin}}-{\alpha}{\epsilon},
\end{equation}\label{def:approx_sol2}}
We know any $(\alpha,\epsilon)-$saddle point, denoted by $(\bar\x,\bar S)$, has the following properties:
%%%%%%%%%%%%%%%%%%%%%%%%%%%%
\begin{enumerate}
    \item $ \underline{\phi}(\bar S)>{\alpha}.{\rm{OPT}}_{maxmin}-\epsilon$
    \item $\bar{\phi}(\bar\x)<\frac{1}{\alpha} .{\rm{OPT}}_{minmax}+\frac{\epsilon}{\alpha}$
\end{enumerate}
%%%%%%%%%%%%%%%%%%%%%%%%%%%%
This is due to the fact that we have:
\begin{enumerate}
    \item $\min\limits_{\x\in \mathcal{X}}f(\x,\bar S)\leq \min\limits_{\x\in \mathcal{X}}\max\limits_{S\in \mathcal{
    I}}f(\x,S)={\rm{OPT}}_{minmax}$
    \item ${\rm{OPT}}_{maxmin}=\max\limits_{S\in \mathcal{
    I}}\min\limits_{\x\in \mathcal{X}}f(\x,S)\leq \max\limits_{S\in \mathcal{
    I}} f(\bar\x,S)$
\end{enumerate}

these two conditions imply that by finding an $(\alpha,\epsilon)-$saddle point we find an $\alpha-$approximate solution for the minimax problem \eqref{eq:main_problem} and a $\frac{1}{\alpha}-$approximate solution for the max-min problem \eqref{eq:main_problem22}. In order to prove finding $(\alpha,\epsilon)-$ saddle point is NP-hard, we will prove that finding approximate solution for maxmin convex-submodular is NP-hard. We do this establishing a connection between this problem and the problem of robust submodular maximization through following result stated and proved in \citep{krause2008robust}.

Consider monotone-submodular functions $f_1, f_2, \dots f_n$ and the following robust  submodular maximization problem:
\begin{equation}
  {\rm{OPT}_1}=\max_{|S|\leq k}\min_{i\in [n]} f_i(S). 
\end{equation}
Solving this problem up to approximation factor is NP-hard, i.e. finding a solution  $S$ such that $\max_i f_i(S)\geq\alpha \text{OPT}_1$  is an NP-hard task for any $\alpha > 0$.

Now, consider the following problem: 
\begin{equation}\label{eq:controbustsubmodular}
   {\rm{OPT}_2}=\max_{|S|\leq k}\,
   \min_{\substack{x\in{\mathbbm{R}}^n,x\geq 0\\x^{T}\mathbbm{1}= 1}}\sum_{i=1}^{n}
   x_i.f_i(S) 
\end{equation}
where $\mathbbm{1}$ is vector of all ones and $x=[x_1,x_2,x_3\dots x_n]^T$.
For this problem, it is easy to verify that $ {\rm{OPT}_1}= {\rm{OPT}_2}$ since for every set $S\in V$ we have  $\min_{i\in [n]} f_i(S)=\min\limits_{\substack{x\in{\mathbbm{R}}^n,x\geq 0\\x^{T}\mathbbm{1}= 1}}\sum_{i=1}^{n} x_i.f_i(S)$. Therefore, finding a  $\alpha-$approximate solution for problem in \eqref{eq:controbustsubmodular} is NP-hard. Problem \eqref{eq:controbustsubmodular} is max-min convex-submodular optimization which means max-min convex-submodular optimization is NP-hard in general. We show that by finding $(\alpha,\epsilon)-$saddle point we can provide $\frac{1}{\alpha}-$approximate solution for max-min problem; therefore, since we proved finding $\frac{1}{\alpha}-$approximate solution for max-min problem is NP-hard, finding $(\alpha,\epsilon)-$saddle point is NP-hard too.
%%%%%%%%%%%%%%%%%%%%%%%%%%%%
\subsection{Proof of Theorem \ref{thm:approx_neg}}
%%%%%%%%%%%%%%%%%%%%%%%%%%%%%
{

{
Before stating this proof, let us explain what we mean by ``NP-hard'' for the considered setting. We note that an algorithm for Problem \eqref{eq:main_problem} is supposed to search for an approximate solution only in  $\mathcal{X}$ (i.e., in terms of the variable $\x$), and for this, it will require some information about the values $f(\x,S)$. However, for every fixed $\x$, there may be restrictions on obtaining some specific values of $f(\x,S)$. For example, finding the exact value of $\bar{\phi}(\x)$ can in general be NP-hard (as maximizing a monotone-submodular function beyond $(1-1/e)$ approximation is hard). In order to appropriately address these restrictions, we will view our setting as a procedure between the algorithm and an oracle that we now describe below.

Upon receiving an input point $\x_{in}\in \mathcal{X}$, the oracle chooses based on this input a set $S_{out}$ such that $|S_{out}|\leq k$, and returns all the following information: the set $S_{out}$, the value $f(\x_{in},S_{out})$, and the gradient of $f(\x_{in},S_{out})$ with respect to $\x$ at the point $(\x_{in},S_{out})$. The only restriction on the oracle is that it is a polynomial-time oracle, i.e. the oracle's procedure to find the set $S_{out}$ requires poly-time complexity in terms of the size of the ground-set $|V|$. More precisely, there exists an integer $q>0$ such that for any ground set $V$, the oracle uses at most $\mathcal{O}(|V|^{q})$ operations to find the output set $S_{out}$ corresponding to an input $\x_{in}$. Note that we do not put any restriction on what the oracle does apart from having poly-time complexity; e.g. it could output a greedy solution, or it could output a random set $S_{out}$, or it could do any other procedure. We call such an oracle a polynomial-time oracle.
Given this choice of the oracle, the algorithm proceeds in $m$ rounds, and in each round $r\in [m]$, it chooses an input point $\x_r\in \mathcal{X}$
 to query from the oracle. Importantly, we consider algorithms that require a polynomial number of rounds in terms of the size of the ground set $V$. More precisely, for any ground set $V$ and $\alpha, \epsilon$, the number of rounds of the algorithm is at most $c(\alpha,\epsilon)|V|^{q}$
 where $q>0$ is an absolute constant and $c(\alpha,\epsilon)$ is another constant that only depends on $\alpha$ and $\epsilon$. We call such an algorithm a \textbf{polynomial-time algorithm}. Next, we show that no polynomial-time algorithm is capable of finding an $(\alpha,\epsilon)$-approximate of  \eqref{eq:main_problem} for $\alpha>1-1/e$.
 
}

In the following, we assume for simplicity that $\epsilon = 0$. The proof can be trivially extended to any value of $\epsilon$, as we will explain at the end of the proof. Recall from the statement of the theorem that $\alpha = 1-1/e + \gamma$ for a fixed constant $\gamma > 0$. 

 We know for a fact that monotone-submodular maximization beyond the $(1-1/e)$-approximation in NP-hard  \citep{krause2014submodular}. I.e. unless P = NP, for any integer $q>0$ there exists a monotone submodular function $g_1: 2^V \to \mathbb{R}_+$ and an integer $k \leq |V|$ such that finding a set $S$ with cardinally $k$ where $g_1(S) \geq (1-1/e + \gamma/3) \max_{|S| \leq k} g_1(S)$ requires computing more than $|V|^q$ function values (i.e. complexity is larger than $|V|^q$). Consider such a function $g_1$ and the choice of $k$, and define $\text{OPT}_{g_1} = \max_{|S| \leq k} g_1(S) $. We also define another function $g_2: 2^V \to \mathbb{R}_+$ as follows:  $g_2(S) =  \min\{g_1(S), (1-1/e + \gamma/3) \text{OPT}_{g_1}\}$. It is important to note that   finding a set $|S| \leq k$ such that $g_1(S) \neq g_2(S)$ requires complexity larger than $|V|^q$ (also note that the choice of $q$ is arbitrary here, i.e. for every $q$ there exists a $g_1$, etc.).   
 
Consider the integer $n = \lceil 4/\gamma + 1 \rceil$ and let  $\mathcal{X}$ be the $n$-dimensional simplex, i.e.   $$\mathcal{X} = \{\mathbf{x}=({x}_1, \cdots, {x}_n) \text{ s.t. }  \sum_{i=1}^n x_i = 1 \, \& \, x_i \geq 0 \,\, \forall{i = 1, \cdots, n} \}.$$ 
For $j \in \{1, \cdots, n\}$ we define $f_j(\mathbf{x},S): \mathcal{X} \times 2^V \to \mathbb{R}_+$ as 

$$f_j(\mathbf{x},S) = \sum_{i\in [n], i \neq j} {x}_i g_1(S) + {x}_j g_2(S)$$

We note a few facts about each of the functions $f_j(\mathbf{x},S)$: 

(i) Any $(\alpha, 0)$ approximate solution for the function $f_j$ has the property that ${x}_j > 1/n + \gamma/4$. This is simply because for any $(\alpha, 0)$-approximate solution $\mathbf{x} = ({x}_1, \cdots, {x}_n)$ we have
$ \alpha \bar{\phi}(\mathbf{x}) \leq \min_{\mathbf{x} \in \mathcal{X}} \max_{|S| \leq k} f_j(\mathbf{x},S)$, and hence

$$ \alpha  \left({x}_j (1-1/e + \gamma/3) + (1-{x}_j) \right) {\rm{OPT}}_{g_1} \leq (1-1/e + \gamma/3) {\rm{OPT}}_{g_1} $$

From the above inequality (and by noting that $\gamma \in (0,1]$) we can always deduce that ${x}_j > \gamma/2$, and thus ${x}_j > 1/n + \gamma/4$. 

(ii) Given a polynomial-time oracle, we can not distinguish between the functions $f_1,\dots,f_n$ using a query from the oracle. This is because the oracle can not find a set $S$ with carnality at most $k$ for which $g_1(S) \neq g_2 (S)$ (as finding that set by the oracle is intractable), and thus, for the set $S_{\rm out}$ that the oracle finds, the outcome of the oracle will be the function value $f(\mathbf{x}_{\rm in}, S_{\rm out}) = g_1(S_{\rm out})$ and $\nabla_{\mathbf{x}} f(\mathbf{x}_{\rm in}, S_{\rm out}) = g_1(S_{\rm out}) \mathbf{1}_n$ where $\mathbf{1}_n$ is the all-ones vector of dimension $n$. These outputs bear absolutely no information about the index $j$. 
 
 Given the above facts, we are now ready to finalize the proof. Consider the scenario where the index $j$ is chosen uniformly at random inside the set $[n]$, and the algorithm aims at finding an approximate solution of the function $f_j$. Note that the choice of $j$ is hidden to the algorithm. Now, given fact (ii) above, if both the algorithm and oracle are polynomial-time, then in all the rounds and queries, there will be absolutely no information revealed about the index $j$. As a result, the mutual information between the outcome of the queries and the index $j$ will be zero. 
 
 On the other hand, from fact (i) above, if an algorithm can find an $(\alpha, 0)$-approximate solution, we claim that the solution is informative about the index $j$. More precisely, given the solution that the algorithm has found, we can infer the hidden index $j$ using the following procedure: the algorithm's solution $\mathbf{x} = (x_1,\cdots, x_n)$ can be viewed as a probability distribution over the set $\{1,\cdots, n\}$. As a result, if we use this probability distribution  to  draw an integer $\hat{j}$ from the set $\{1,\cdots, n\}$, then we have $\text{Pr}\{\hat{j} = j\} = x_j >= 1/n + \gamma/4$. Thus, we can decode the index $j$ with a probability that is strictly larger than a random guess. This means that the mutual information of the solution found by the algorithm and the index $j$ is strictly lower-bounded by a positive constant (which only depends on $\gamma$). This contradicts the result of the previous paragraph.

Note that in the above we have assumed that $\epsilon = 0$. For general $\epsilon$, we note that we can always choose the function $g_1$ such that $\text{OPT}_{g_1}$ is sufficiently large. As a result, we can write $\epsilon = \epsilon' \times \text{OPT}_{g_1}$ where $\epsilon'$ can be made arbitrarily small. Hence, proving hardness for obtaining an $(\alpha, \epsilon)$-approximate becomes equivalent to proving harness for obtaining an $(\alpha/(1+\epsilon'), 0)$ approximate solution.  The conclusion is now immediate since out proof above works for any $\alpha = 1-1/e + \gamma$ and $\epsilon'$ can be made arbitrarily small by making $\text{OPT}_{g_1}$ sufficiently large.
}
 \subsection{Proof of Theorem~\ref{thm:GG2}: Gradient Greedy(GG) Convergence}\label{sec:Gradient Greedy convergence}
 %%%%%%%%%%%%%%%%%%%%%%%%%%%%
 \proof
 From Assumption~\ref{assum:smooth}, we know $\nabla f(\x,S)$ for all $\x,S$ is $L$-Lipschitz. which results in  $\nabla^2 f(\x,S)	\preceq L\mathbf{I}$, where $\mathbf{I}$ is identity matrix. We can thus write for any $\mathbf{y}$:
 \begin{align}\label{eq:gg1}
     f(\mathbf{y},S)\leq f(\x,S) + \langle\nabla f(\x,S), \mathbf{y}-\x \rangle + \frac{L}{2}\|\mathbf{y}-\x\|^2 
 \end{align}
 Now, define $\x^{+}_S=\x-\alpha \nabla f(\x,S)$. By substituting it in \eqref{eq:gg1} we obtain:
 \begin{align}
     f(\x^{+}_S,S)\leq f(\x,S)+\nabla f(\x,S)^T(\x^{+}_S-\x)+\frac{L}{2}\|\x^{+}_S-\x\|^2
     \\=f(\x,S)-\nabla f(\x,S)^T(\alpha \nabla f(\x,S))+\frac{L}{2}\|-\alpha \nabla f(\x,S)\|^2
     \\=f(\x,S)-(\alpha-\frac{L\alpha^2}{2})\|\nabla f(\x,S)\|^2
 \end{align}
 if $\alpha\leq \frac{1}{L}$ we know that $-(\alpha-\frac{L\alpha^2}{2})\leq -\frac{\alpha}{2}$
 which means:
  \begin{align}\label{eq:gg2}
     f(\x^{+}_S,S)\leq f(\x,S)-\frac{\alpha}{2}\|\nabla f(\x,S)\|^2
 \end{align}
Consequently,  in each step of gradient descent objective value decreases.
 Then, for every point $\tilde{\x}$ we can write using convexity:
 \begin{align}\label{eq:gradg2.5}
     f(\x,S)\leq f(\tilde{\x},S)+ \nabla f(\x,S)^T(\x-\tilde{\x})
 \end{align}
 Combining \eqref{eq:gg2} and \eqref{eq:gradg2.5} we have:
 
 \begin{align}
     f(\x^{+}_S,S)-f(\tilde{\x},S)&\leq  \nabla f(\x,S)^T(\x-\tilde{\x})-\frac{\alpha}{2}\|\nabla f(\x,S)\|^2
     \\&=\frac{1}{2\alpha}(2\alpha \nabla f(\x,S)^T(\x-\tilde{\x})-\alpha^2\|\nabla f(\x,S)\|^2\notag\\&\quad-\|\x-\tilde{\x}\|^2+\|\x-\tilde{\x}\|^2)
     \\&=\frac{1}{2\alpha}(-\|\x-\alpha\nabla f(\x,S)-\tilde{\x}\|^2+\|\x-\tilde{\x}\|^2)
     \\&=\frac{1}{2\alpha}(-\|\x^{+}_S-\tilde{\x}\|^2+\|\x-\tilde{\x}\|^2)\label{eq:gg3}
 \end{align}
 Now, by using \eqref{eq:gg3}, if we let $\gamma_t=\alpha$ we have the following inequality for $(\x_t,S_{t})$:
 
 \begin{align}
     f(\x_{t},S_t)-f(\tilde{\x},S_t)\leq \frac{1}{2\alpha}(-\|\x_{t}-\tilde{\x}\|^2+\|\x_{t-1}-\tilde{\x}\|^2)
 \end{align}
 summing up over $t$ we have:
 \begin{align}
    \sum_{t=1}^{T} f(\x_{t},S_t)-f(\tilde{\x},S_t)\leq \frac{1}{2\alpha}(\|\x_{0}-\tilde{\x}\|^2)
 \end{align}
 
 our set of continuous variable is bounded which means $\|\x\|^2\leq H$; this results:
 
 \begin{align}\label{eq:gg6}
    \sum_{t=1}^{T} f(\x_{t},S_t)-f(\tilde{\x},S_t)\leq \frac{H}{2\alpha}
 \end{align}
 
 Also, from greedy update we have for every $S$(check \citep{krause2014submodular}):
 \begin{align}\label{eq:gg4}
 f(\x_{t-1},S)-\frac{f(\x_{t-1},S_{t})}{1-\frac{1}{e}}\leq 0
 \end{align}
 Now, using the Lipschitz condition (consequence of Assumption \ref{assum:bounded_grad}):
 \begin{align}\label{eq:gg5}
 |f(\x_{t},S_t)-f(\x_{t-1},S_t)|\leq M\|\x_t-\x_{t-1}\|\leq M\alpha\|\nabla f(\x_{t-1},S_t)\|\leq M^2\alpha
 \end{align}
 Putting \eqref{eq:gg4} and \eqref{eq:gg5} together:
 \begin{align}
 ({1-\frac{1}{e}})f(\x_{t},S)-{f(\x_t,S_{t})}\leq 2M^{2}\alpha
 \end{align}
 and summing over $t$ we have:
 \begin{align}\label{eq:gg7}
\sum_{t=1}^{T} ({1-\frac{1}{e}})f(\x_{t},S)-{f(\x_t,S_{t})}\leq 2M^{2}\alpha T
 \end{align}
 From \eqref{eq:gg6} and \eqref{eq:gg7} we can then obtain the following:

 \begin{align}\label{eq:gg81}
\sum_{t=1}^{T} ({1-\frac{1}{e}})f(\x_{t},S)-f(\tilde{\x},S_t)\leq 2M^{2}\alpha T+\frac{H}{2\alpha}
 \end{align}
and finally:
  \begin{align}\label{eq:gg8}
\frac{\sum_{t=1}^{T} \alpha (({1-\frac{1}{e}})f(\x_{t},S)-f(\tilde{\x},S_t))}{ \sum_{t=1}^{T}\alpha}\leq \frac{2M^{2}\alpha T+\frac{H}{2\alpha}}{T}
 \end{align}
  %%%%%%%%%%%%%%%%%%
    
  %%%%%%%%%%%%%%%%%%
From convexity we have:
 \begin{equation}
        f(\frac{1}{T}\sum_{t=1}^{T}\x_{t},S)  \leq\frac{\sum_{t=1}^{T}f(\x_{t-1},S)}{T}
 \end{equation}
 which results in:
  \begin{align}\label{eq:gg9}
({1-\frac{1}{e}}) f(\frac{1}{T}\sum_{t=1}^{T}\x_{t},S) -\frac{\sum_{t=1}^{T} \alpha (f(\tilde{\x},S_t))}{ \sum_{t=1}^{T}\alpha}\leq \frac{2M^{2}\alpha T+\frac{H}{2\alpha}}{T}
 \end{align}
 For $\x^{*}$,
 we know that
 $\min_\x\max_S f(\x,S)=\max f(\x^*,S)\geq f(\x^*,S_t)$.
 Now in \eqref{eq:gg9} we let $\tilde{\x}=\x^*$ and write:
 \begin{align}\label{eq:gg10}
({1-\frac{1}{e}})\max_S f(\frac{1}{T}\sum_{t=1}^{T}\x_{t},S) -\min_x\max_S f(\x,S)\leq \frac{2M^{2}\alpha T+\frac{H}{2\alpha}}{T}
 \end{align}
 Letting  $\alpha=\frac{1}{\sqrt{T}}$ we obtain:
 \begin{align}\label{eq:gg11}
({1-\frac{1}{e}})\max_S f(\frac{1}{T}\sum_{t=1}^{T}\x_{t},S) -\min_\x\max_S f(\x,S)\leq \frac{2M^{2}+\frac{H}{2}}{\sqrt{T}}
 \end{align}
 Finally, if we define $K=2M^{2}+\frac{H}{2}$ and let $T=\frac{K^2}{\epsilon^2}$; then $\x_{sol}=\frac{1}{T}\sum_{t=1}^{T}\x_{t}$ is a $(1-1/e,\epsilon)$- approximate minimax solution.
  %%%%%%%%%%%%%%%%%%
 
\subsection{Proof of Theorem~\ref{thm:GRG2}:  Gradient Replacement-greedy (GRG) Convergence}\label{sec:Gradient Replacement-greedy convergence}

\proof
Let $g$ be a monotone-submodular function, and consider sets $B,S\subseteq V$ with size $k$. Define  $e^{*}=\argmax_{e\in S} g(S\setminus e)-g(S)$, and $v^{*}=\argmax_{v\in V} g(S\cup v\setminus e^{*})-g(S\cup v\setminus e^{*})$.  We have:
\begin{align}\label{eq:rg1proof}
g(S\cup v^{*}\setminus e^{*})-g(S)\geq \frac{1}{k}\sum_{v\in B}g(S\cup v\setminus e^{*})-g(S)
\notag\\= \frac{1}{k}\sum_{v\in B}(g(S\cup v\setminus e^{*})-g(S\cup v)+g(S\cup v)-g(S))
\end{align}
where the first inequality comes from the definition of $v^{*}$. We know that for a monotone-submodular function $g$ we have
$g(B\cup S)-g(S)\leq \sum_{v\in B}(g(S\cup v)-g(S))$ for any choice of $B,S$ \citep{stan2017probabilistic}; which results in:

\begin{align}\label{eq:rg2proof}
\frac{1}{k}\sum_{v\in B}(g(S\cup v)-g(S))
\geq \frac{1}{k}(g(B\cup S)-g(S))\geq \frac{1}{k}(g(B)-g(S))
\end{align}
Here, the first inequality is due to  submodularity and the second inequalities is due to monotonicity. Also, we have:
\begin{align}\label{eq:rg3proof}
    \frac{1}{k}\sum_{v\in B}(g(S\cup v)-g(S\cup v\setminus e^{*}))&\leq \frac{1}{k}\sum_{v\in B}(g(S)-g(S\setminus e^{*}))=g(S)-g(S\setminus e^{*})
    \notag\\&\leq\frac{1}{k}\sum_{e\in S}(g(S)-g(S\setminus e))\leq \frac{1}{k}g(S)
\end{align}
where the first and second inequality comes from submodularity. Combining \eqref{eq:rg1proof},\eqref{eq:rg2proof}, and \eqref{eq:rg3proof} we have that for every set $B$ of size $k$:
\begin{align}\label{eq:rg4proof}
g(S\cup v^{*}\setminus e^{*})-g(S)\geq \frac{1}{k}(g(B)-2g(S))
\end{align}

If we apply \eqref{eq:rg4proof} for the replacement greedy update in Gradient Replacement-greedy(GRG) algorithm, we obtain:

\begin{align}\label{eq:GRG1}
    f( \x_t,{S}_t)- f( \x_t, S_{t-1})\geq \frac{1}{k}\left(
    f( \x_t, S)- 2f( \x_t, S_{t-1})
    \right)
\end{align}

and hence
\begin{equation} \label{eq:GRG2}
f( \x_t, S) - 2f(\x_{t},  S_{t}) \leq
(1-\frac{2}{k}) (f( \x_t, S) - 2f( \x_{t}, {S}_{t-1}))
\end{equation}
Note that as $f$ is $M$-Lipschitz we have for every $S$ (consequence of Assumption \ref{assum:bounded_grad}):

\begin{equation} \label{eq:GRG3}
|f( \x_t, S)-f( \x_{t-1}, S)|\leq M\|\x_t-\x_{t-1}\|\leq M\alpha\|\nabla f(\x_{t-1},S_t)\|\leq M^2\alpha
\end{equation}
Combining \eqref{eq:GRG2} and \eqref{eq:GRG3} we obtain that
\begin{equation} \label{eq:GRG2.2}
f( \x_t, S) - 2f(\x_{t},  S_{t}) \leq
(1-\frac{2}{k}) (f( \x_{t-1}, S) - 2f( \x_{t-1}, {S}_{t-1})+3M^{2}\alpha)
\end{equation}
Using a recursive argument we can show that 
\begin{align} \label{eq:GRG2.3}
f( \x_t, S) - 2f(\x_{t},  S_{t}) \leq
(1-\frac{2}{k})^t (f( \x_{0}, S) - 2f( \x_{0}, {S}_{0})) +\sum_{m=1}^{t}(1-\frac{2}{k})^m 3M^{2}\alpha
\end{align}
Now since $f( \x_{0}, {S}_{0})$ is non-negative, we can eliminate $-f( \x_{0}, {S}_{0})$ from the right hand side. Using this observation and by simplifying the geometric sum we obtain that
\begin{align}
f( \x_t, S) - 2f(\x_{t},  S_{t})
\leq (1-\frac{2}{k})^t f( \x_{0}, S)  + 3M^{2}\alpha \frac{k}{2}
\end{align}
Now, note that  $(1-\frac{2}{k})^t$ is bounded above by $e^{-\frac{2t}{k}}$ and therefore we have 
\begin{align}\label{sth}
f( \x_t, S) - 2f(\x_{t},  S_{t})
\leq Ae^{-\frac{2t}{k}}+3M^{2}\alpha \frac{k}{2},
\end{align} 
where $A$ is an upper bound for function value at point zero, $f( 0, S)\leq A$.
Now, from the analysis of gradient descent similar to \eqref{eq:gg6}, we have:
\begin{align}\label{eq:GRG3.2}
    \sum_{t=1}^{T} f(\x_{t},S_t)-f(\tilde{\x},S_t)\leq \frac{H}{2\alpha}
 \end{align}
 
Combining this inequality with \eqref{sth} we have:
 \begin{align}\label{eq:GRG3.3}
    \sum_{t=1}^{T} \frac{1}{2}f( \x_t, S)-f(\tilde{\x},S_t)&\leq \frac{H}{2\alpha}+\frac{\sum_{t=1}^{T}Ae^{-\frac{2t}{k}}}{2}+\frac{3TM^{2}\alpha k}{4}\notag\\&\leq \frac{K}{2\alpha}+\frac{Ae^{-\frac{2}{k}}}{2(1-e^{-\frac{2}{k}})}+\frac{3TM^{2}\alpha k}{4}
 \end{align}
Thus, choosing the parameters $\alpha=\frac{1}{\sqrt{T}}$ will lead to
 \begin{align}\label{eq:GRG4}
   \frac{1}{T} \sum_{t=1}^{T} \frac{1}{2}f( \x_t, S)-f(\tilde{\x},S_t)\leq \frac{H}{2\sqrt{T}}+\frac{Ae^{-\frac{2}{k}}}{2T(1-e^{-\frac{2}{k}})}+\frac{3M^{2} k}{4\sqrt{T}}\leq \frac{K}{\sqrt{T}}
 \end{align}
 where $K$ is some constant.

 In summary, we have obtained the following relation that will be used to drive the guarantee for the minimax problem:
 \begin{align}\label{eq:GRG4final}
   \frac{1}{T} \sum_{t=1}^{T} \frac{1}{2}f( \x_t, S)-f(\tilde{\x},S_t)\leq \frac{K}{\sqrt{T}}
 \end{align}
   %%%%%%%%%%%%%%%%%%
  %%%%%%%%%%%%%%%%%%
  We know because of convexity we have:
 \begin{equation}\label{eq:GRG5}
        f\left(\frac{1}{T}\sum_{t=1}^{T}\x_{t-1},S\right)  \leq\frac{\sum_{t=1}^{T}f(\x_{t-1},S)}{T}
 \end{equation}

 Now combining \eqref{eq:GRG4final} and \eqref{eq:GRG5} we have: 
 \begin{align}\label{eq:GRG61}
 \frac{1}{2}f\left(  \frac{1}{T} \sum_{t=1}^{T}  \x_t, S\right)-  \frac{1}{T} \sum_{t=1}^{T} f(\tilde{\x},S_t)\leq \frac{K}{\sqrt{T}}
 \end{align}
 
  Also for  $\x^{*}$
we have $\min_\x\max_S f(\x,S)=\max_S f(\x^*,S)\geq f(\x^*,S_t)$. By using $\tilde{\x}=\x^*$ we can write:
  \begin{align}\label{eq:GRG6}
 \frac{1}{2}f\left(  \frac{1}{T} \sum_{t=1}^{T}  \x_t, S\right)-  \min_\x\max_S f(\x,S)&\leq \frac{1}{2}f\left(  \frac{1}{T} \sum_{t=1}^{T}  \x_t, S\right)-  \frac{1}{T} \sum_{t=1}^{T} f({\x}^{*},S_t)\notag\\&\leq \frac{K}{\sqrt{T}}
 \end{align}
  Let $T=\frac{K^2}{\epsilon^2}$; then $\x_{sol}=\frac{1}{T}\sum_{t=1}^{T}\x_{t}$ is a $(1/2,\epsilon)$- approximate minimax solution.
  %%%%%%%%%%%%%%%%%%
 
\subsection{Proof of Theorem \ref{thm:EGG2}: Extra-gradient Greedy(EGG) Convergence  }\label{sec:Extra-Gradient Greedy Convergence}
Consider the Extra-gradient Greedy method, we can write the following equations to find the bound on convergence of $\x$:
\begin{align}
&\| \hat{\x} -\x\|^2\nonumber\\
& \leq \| \x_t -\x -\gamma_t \nabla_\x f(\x_t,S_t)\|^2 \nonumber\\
& =  \| \x_t -\x\|^2 -2 \gamma_t \nabla_\x f(\x_t,S_t)^\top (\x_t -\x )+ \|\hat{\x_t}-\x_t \|^2\nonumber\\
& = \| \x_t -\x\|^2 -2 \gamma_t \nabla_\x f(\x_t,S_t)^\top (\hat{\x}_t -\x )+ \|\hat{\x}_t-\x_t \|^2
\notag\\&\quad+2 \gamma_t \nabla_\x f(\x_t,S_t)^\top (\hat{\x}_t -\x_t ) \nonumber\\
& = \| \x_t -\x\|^2 -2 \gamma_t \nabla_\x f(\x_t,S_t)^\top (\hat{\x}_t -\x )+ \|\hat{\x}_t-\x_t \|^2
+2(\x_t-\hat{\x}_t )^\top (\hat{\x}_t -\x_t ) \nonumber\\
& = \| \x_t -\x\|^2 -2 \gamma_t \nabla_\x f(\x_t,S_t)^\top (\hat{\x}_t -\x )+ \|\hat{\x}_t-\x_t \|^2
-2\|\hat{\x}_t -\x_t \|^2 
\end{align}
Hence, we have 
\begin{align}\label{popo}
2\gamma_t \nabla_\x f(\x_t,S_t)^\top (\hat{\x}_t -\x )\leq \| \x_t -\x\|^2-\| \hat{\x}_t -\x\|^2-\|\hat{\x}_t -\x_t \|^2 
\end{align}

Similarly we can show that 
\begin{align}
&
\| \x_{t+1} -\x\|^2\nonumber\\
& \leq \| \x_t -\x -\gamma_t \nabla_\x f(\hat{\x_t},\hat{S}_t)\|^2 \nonumber\\
& =  \| \x_t -\x\|^2 -2 \gamma_t \nabla_\x f(\hat{\x_t},\hat{S}_t)^\top (\x_t -\x )+ \|\x_{t+1}-\x_t \|^2\nonumber\\
& = \| \x_t -\x\|^2 -2 \gamma_t \nabla_\x f(\hat{\x_t},\hat{S}_t)^\top (\x_{t+1} -\x )+ \|\x_{t+1}-\x_t\|^2
\notag\\&\quad+2 \gamma_t \nabla_\x f(\hat{\x_t},\hat{S}_t)^\top (\x_{t+1} -\x_t ) \nonumber\\
& = \| \x_t -\x\|^2 -2 \gamma_t   \nabla_\x f(\hat\x_{t},\hat{S}_t)^\top (\hat{\x}_t -\x )+ \|\x_{t+1}-\x_t \|^2
\notag\\&\quad+2(\x_t-\x_{t+1}  )^\top (\x_{t+1} -\x_t ) \nonumber\\
& = \| \x_t -\x\|^2 -2 \gamma_t \nabla_\x f(\hat{\x_t},\hat{S}_t)^\top (\x_{t+1}-\x )
\notag\\&\quad+ \|\x_{t+1}-\x_t \|^2
-2\|\x_{t+1} -\x_t \|^2 
\end{align}

Hence, we have 
\begin{align}\label{eq:extragradient main}
2\gamma_t \nabla_\x f(\hat{\x_t},\hat{S}_t)^\top (\x_{t+1}-\x )\leq \| \x_t -\x\|^2-\| \x_{t+1} -\x\|^2-\|\x_{t+1} -\x_t \|^2 
\end{align}

Now note that we can write $2\gamma_t \nabla_\x f(\hat{\x}_t,\hat{S}_t)^\top (\hat{\x}_t -\x)$ as 
\begin{align}\label{sssddd}
 &2\gamma_t \nabla_\x f(\hat{\x}_t,\hat{S}_t)^\top (\hat{\x}_t -\x)\\
&=\  2\gamma_t \nabla_\x f(\hat{\x}_t,\hat{S}_t)^\top (\hat{\x}_t -\x_{t+1} )+ 2\gamma_t \nabla_\x f(\hat{\x}_t,\hat{S}_t)^\top (\x_{t+1} -\x )  \\
&=\ 2\gamma_t \nabla_\x f(\hat{\x}_t,\hat{S}_t)^\top (\hat{\x}_t -\x_{t+1} )+ 2\gamma_t \nabla_\x f(\hat{\x}_t,\hat{S}_t)^\top (\x_{t+1} -\x )  \\
&\qquad + 2\gamma_t \nabla_\x f({\x}_t,{S}_t)^\top (\hat{\x}_t-\x_{t+1} )- 2\gamma_t \nabla_\x f({\x}_t,{S}_t)^\top (\hat{\x}_t-\x_{t+1})\\
&=\  2\gamma_t \nabla_\x f({\x}_t,{S}_t)^\top (\hat{\x}_t-\x_{t+1} ) + 2\gamma_t \nabla_\x f(\hat{\x}_t,\hat{S}_t)^\top (\x_{t+1}-\x )\\
&\qquad +2\gamma_t \left(\nabla_\x f(\hat{\x}_t,\hat{S}_t)- \nabla_\x f({\x}_t,{S}_t)\right)^\top (\hat{\x}_t-\x_{t+1})\\
&\leq\  
\| \x_t -\x_{t+1}\|^2-\| \hat{\x_t} -\x_{t+1}\|^2-\|\hat{\x}_t -\x_t \|^2 \\ 
&\qquad +  \| \x_t -\x\|^2-\| \x_{t+1} -\x\|^2-\|\x_{t+1} -\x_t \|^2 \\
&\qquad +2\gamma_t \left(\nabla_\x f(\hat{\x}_t,\hat{S}_t)- \nabla_\x f({\x}_t,{S}_t)^\top (\hat{\x}_t-\x_{t+1})\right)\\
&=\  
-\| \hat{\x_t} -\x_{t+1}\|^2-\|\hat{\x}_t -\x_t \|^2  +  \| \x_t -\x\|^2-\| \x_{t+1} -\x\|^2 \\
&\qquad +2\gamma_t \left(\nabla_\x f(\hat{\x}_t,\hat{S}_t)- \nabla_\x f({\x}_t,{S}_t)\right)^\top (\hat{\x}_t-\x_{t+1}),
\end{align}
where the inequality follows from the results in \eqref{popo} and \eqref{eq:extragradient main}. 

Next we derive an upper bound for the inner product \newline$\left(\nabla_\x f(\hat{\x}_t,\hat{S}_t)- \nabla_\x f({\x}_t,{S}_t)\right)^\top (\hat{\x}_t-\x_{t+1})$ using the smoothness of the function $f$, i.e., 
\begin{align*}
&\left(\nabla_\x f(\hat{\x}_t,\hat{S}_t)- \nabla_\x f({\x}_t,{S}_t)\right)^\top (\hat{\x}_t-\x_{t+1})\\
&\leq \| \nabla_\x f(\hat{\x}_t,\hat{S}_t)- \nabla_\x f({\x}_t,{S}_t)\| \|\hat{\x}_t-\x_{t+1}\| \\
&\leq \left(\| \nabla_\x f(\hat{\x}_t,\hat{S}_t)- \nabla_\x f(\hat{\x}_t,{S}_t)\|+\| \nabla_\x f(\hat{\x}_t,{S}_t)- \nabla_\x f({\x}_t,{S}_t)\|\right) \|\hat{\x}_t-\x_{t+1}\| \\
&\leq \left(L_{\x,S}\| \hat{S}_t-{S}_t\|+L_{\x,\x}\| \hat{\x}_t- {\x}_t\|\right) \|\hat{\x}_t-\x_{t+1}\| \\
\end{align*}
Now to complete our upper bound we need to bound $\| \hat{S}_t-{S}_t\|$ which can be done as
\begin{align}
\| \hat{S}_t-{S}_t\| \leq \phi \|\hat{\x}_t-\x_t\|+\sigma
\end{align}
The above relation holds because  for every two feasible set we have $\| A-B\|\leq 2k$
;therefore, if we let $\sigma=2k$ and $\phi=1$ the above condition is always true. Considering this result we obtain that 
\begin{align*}
&\left(\nabla_\x f(\hat{\x}_t,\hat{S}_t)- \nabla_\x f({\x}_t,{S}_t)\right)^\top (\hat{\x}_t-\x_{t+1})\notag\\&\quad\leq \left(L_{\x,S}\phi +L_{\x,\x}\right) \|\hat{\x}_t-\x_t\|\|\hat{\x}_t-\x_{t+1}\| 
+L_{\x,S}\sigma \|\hat{\x}_t-\x_{t+1}\| 
\end{align*}
Applying this upper bound into \eqref{sssddd} implies that
\begin{align*}
 &2\gamma_t \nabla_\x f(\hat{\x}_t,\hat{S}_t)^\top (\hat{\x}_t -\x)\\
&\leq\  -\| \hat{\x_t} -\x_{t+1}\|^2-\|\hat{\x}_t -\x_t \|^2  +  \| \x_t -\x\|^2-\| \x_{t+1} -\x\|^2\\ &\qquad +2\gamma_t \left(L_{\x,S}\phi +L_{\x,\x}\right) \|\hat{\x}_t-\x_t\|\|\hat{\x}_t-\x_{t+1}\| 
+2 \gamma_t L_{\x,S}\sigma \|\hat{\x}_t-\x_{t+1}\| \\
&\leq \| \x_t -\x\|^2-\| \x_{t+1} -\x\|^2\\ 
&\qquad +[  -\| \hat{\x_t} -\x_{t+1}\|^2-\|\hat{\x}_t -\x_t \|^2 \\&\qquad+2\gamma_t \left(L_{\x,S}\phi +L_{\x,\x}\right) \|\hat{\x}_t-\x_t\|\|\hat{\x}_t-\x_{t+1}\|] +2 \gamma_t L_{\x,S}\sigma \|\hat{\x}_t-\x_{t+1}\| \\
&\leq \| \x_t -\x\|^2-\| \x_{t+1} -\x\|^2  -\| \hat{\x_t} -\x_{t+1}\|^2-\|\hat{\x}_t -\x_t \|^2
\\ 
&\qquad +\gamma_t \left(L_{\x,S}\phi +L_{\x,\x}\right) \|\hat{\x}_t-\x_t\|^2+\gamma_t \left(L_{\x,S}\phi +L_{\x,\x}\right)\|\hat{\x}_t-\x_{t+1}\|^2 \\&\qquad+ 4\gamma_t^2L_{\x,S}^2 \sigma^2 +\frac{1}{4} \|\hat{\x}_t-\x_{t+1}\|^2 \\
&\leq \| \x_t -\x\|^2-\| \x_{t+1} -\x\|^2+ 4\gamma_t^2 L_{\x,S}^2\sigma^2 
\end{align*}
where the third inequality holds because of the fact that $2ab\leq a^2+b^2$, and the last inequality holds since we assume $\gamma_t (L_{\x,S}\phi +L_{\x,\x})\leq 3/4$.

Using this result we have that 
\begin{align*}
2\gamma_t \nabla_\x f(\hat{\x}_t,\hat{S}_t)^\top (\hat{\x}_t -\x)\leq \| \x_t -\x\|^2-\| \x_{t+1} -\x\|^2+4\gamma_t^2 L_{\x,S}^2\sigma^2
\end{align*}
Now by convexity of $f$ with respect to $\x$ we have 
\begin{align*}
\nabla_\x f(\hat{\x}_t,\hat{S}_t)^\top (\hat{\x}_t -\x)\geq f(\hat{\x}_t,\hat{S}_t)-f(\x,\hat{S}_t)
\end{align*}
and therefore 
\begin{align*}
f(\hat{\x}_t,\hat{S}_t)-f(\x,\hat{S}_t)\leq\frac{1}{2\gamma_t}\left( \| \x_t -\x\|^2-\| \x_{t+1} -\x\|^2\right)
+2\gamma_t L_{\x,S}^2\sigma^2
\end{align*}
Moreover we know that 
\begin{align*}
f(\hat{\x}_t,{S}) -  \frac{1}{1-1/e}f(\hat{\x}_t,\hat{S}_t)\leq 0
\end{align*}
Hence,
\begin{align}\label{eq:EGGmain}
f(\hat{\x}_t,\hat{S}_t)-f(\x,\hat{S}_t) &+({1-1/e}) f(\hat{\x}_t,{S}) - f(\hat{\x}_t,\hat{S}_t)\notag\\&\leq \frac{1}{2\gamma_t} \left(\| \x_t -\x\|^2-\| \x_{t+1} -\x\|^2\right) +2\gamma_t L_{\x,S}^2\sigma^2
\end{align}
Let $\gamma_t=\frac{1}{\sqrt{T}}$. Then, since $\|\x\|^2\leq H$, we have
\begin{align}\label{eq:EGGmain2}
\frac{1}{T}\sum_{t=1}^{T}-f(\x,\hat{S}_t) +({1-1/e}) f(\hat{\x}_t,{S})\leq \frac{H}{2\sqrt{T}}  +\frac{2 L_{\x,S}^2\sigma^2}{\sqrt{T}}
\end{align}
Let $K=2 L_{\x,S}^2\sigma^2+K/2$ then we have
\begin{align}\label{eq:EGGmain3}
\frac{1}{T}\sum_{t=1}^{T}-f(\x,\hat{S}_t) +({1-1/e}) f(\hat{\x}_t,{S})\leq \frac{K}{\sqrt{T}}
\end{align}
  %%%%%%%%%%%%%%%%%%
  We know because of convexity
 \begin{equation}\label{eq:EGG5}
        f(\frac{1}{T}\sum_{t=1}^{T}\hat \x_{t-1},S)  \leq\frac{\sum_{t=1}^{T}f(\hat \x_{t-1},S)}{T}
 \end{equation}

 Now combining \eqref{eq:EGGmain3} and \eqref{eq:EGG5} we have
 \begin{align}\label{eq:EGG61}
 (1-1/e)f(  \frac{1}{T} \sum_{t=1}^{T}  \hat \x_t, S)-  \frac{1}{T} \sum_{t=1}^{T} f(\tilde{\x},\hat S_t)\leq \frac{K}{\sqrt{T}}
 \end{align}
 
  Also for  $\x^{*}$,
we have $\min_\x\max_S f(\x,S)=\max_S f(\x^*,S)\geq f(\x^*,S_t)$. We let $\tilde{\x}=\x^*$ and  write
  \begin{align}\label{eq:EGG6}
 (1-1/e)f(  \frac{1}{T} \sum_{t=1}^{T}  \hat \x_t, S)&-  \min_\x\max_S f(\x,S)\\&\leq (1-1/e)f(  \frac{1}{T} \sum_{t=1}^{T}  \hat \x_t, S)-  \frac{1}{T} \sum_{t=1}^{T} f({\x}^{*}, \hat S_t)\leq \frac{K}{\sqrt{T}}
 \end{align}
  Let $T=\frac{K^2}{\epsilon^2}$; then, $\x_{sol}=\frac{1}{T}\sum_{t=1}^{T} \hat \x_{t}$ is an $((1-1/e),\epsilon)$ approximate minimax solution.
  %%%%%%%%%%%%%%%%%%
\subsection{Extra-gradient Replacement-greedy(EGRG) Convergence}\label{sec:Extra-Gradient ReplacementGreedy Convergence}

For the analysis with respect to $\x$, we can show that
\begin{align*}
f(\hat{\x}_t,\hat{S}_t)-f(\x,\hat{S}_t)\leq\frac{1}{2\gamma_t}\left( \| \x_t -\x\|^2-\| \x_{t+1} -\x\|^2\right)
+2\gamma_t L_{\x,S}^2\sigma^2
\end{align*}
therefore:
%%%%%%%%%%%%%%%%%%

  \begin{align}\label{eq:x EgRg}
     \frac{\sum_{t=1}^{T}\gamma_t\left [f(\hat{\x}_t,\hat{S}_t)-f(\x,\hat{S}_t)\right]}{\sum_{t=1}^{T}\gamma_t}\leq\frac{\|\x-{\x}_1\|^2}{2\sum_{t=1}^{T}\gamma_t}+\frac{\sum_{t=1}^{T}2\gamma_t^2 L_{x,S}^2\sigma^2}{\sum_{t=1}^{T}\gamma_t}\leq \frac{K_1}{\sqrt{T}}
    \end{align}
 %%%%%%%%%%%%%%%   

It remains to derive an upper bound for
$f(\hat{\x}_t,{S}) -  2 f(\hat{\x}_t,\hat{S}_t)
$ . According to the update of replacement-greedy method, we can write the following inequalities:

\begin{align}\label{first}
    f(\x_t,\hat{S}_t)- f(\x_t,S_t)\geq \frac{1}{k}\left(
    f(\x_t, S)- 2f(\x_t,S_t)
    \right)
\end{align}
and 
\begin{align}
    f(\hat{\x}_t,S_{t+1})- f(\hat{\x}_t,\hat{S}_t)\geq \frac{1}{k}\left(
    f(\hat{\x}_t, S)- 2f(\hat{\x}_t,\hat{S}_t)
    \right)
\end{align}

Using the second expression, we can write 
\begin{align}
 \frac{1}{k}\left(
    f(\hat{\x}_{t-1}, S)- 2f(\hat{\x}_{t-1},\hat{S}_{t-1})\right)
    \leq f(\hat{\x}_{t-1},S_{t})- f(\hat{\x}_{t-1},\hat{S}_{t-1})
\end{align}
Let $\bar{\phi}(\x_t)  = \max_{|S| \leq k } f(\x_t, S) $; if we assume for every $\x,S$, $\|\nabla_{\x} f(\x,S)\|\leq G$ then we have: 
%%%%%%%%%%%%%%%%
\begin{equation} \label{assumption_lp_phi_bar} 
| \bar{\phi}(\x) - \bar{\phi}(\y)| \leq  G || \x - \y||
\end{equation}
%%%%%%%%%%%%%%%

 hence
\begin{equation} \label{1-2/k}
\bar{\phi}(\hat \x_{t-1}) - 2f(\hat \x_{t-1},  S_{t}) \leq
(1-\frac{2}{k}) (\bar{\phi}(\hat \x_{t-1}) - 2f(\hat \x_{t-1}, \hat{S}_{t-1}))
\end{equation}
Note that 
{\begin{align} \nonumber
    f(\hat \x_{t-1}, {S}_{t}) &\leq  f( \x_{t}, {S}_{t}) + L_x ||\x_{t} - \hat \x_{t-1}|| \\
    &\leq  f( \x_{t}, \hat{S}_{t}) + \gamma_t G^2\\
    &\leq  f( \hat \x_{t}, \hat{S}_{t}) + 2\gamma_t G^2, \label{eq:s and s_hat}
\end{align}}
therefore,
\begin{equation} \label{eq:1-2/k}
\bar{\phi}(\hat \x_{t-1}) - 2f(\hat \x_{t},  \hat S_{t}) \leq
(1-\frac{2}{k}) (\bar{\phi}(\hat \x_{t-1}) - 2f(\hat \x_{t-1}, \hat{S}_{t-1}))+ 4\gamma_t G^2
\end{equation}
Also, note that
\begin{equation} \label{eq:x_t phi}
    | \bar{\phi}(\hat \x_{t}) - \bar{\phi}(\hat \x_{t-1}) |
    \leq  G^2 \gamma_t.
\end{equation}

Putting \eqref{eq:1-2/k}, \eqref{eq:s and s_hat} and \eqref{eq:x_t phi} together, we obtain:
\begin{equation} \label{eq:1-2/k 2}
\bar{\phi}(\hat \x_{t}) - 2f(\hat \x_{t},  \hat S_{t}) \leq
(1-\frac{2}{k}) (\bar{\phi}(\hat \x_{t-1}) - 2f(\hat \x_{t-1}, \hat{S}_{t-1}))+ 5\gamma_t G^2
\end{equation}

 let $\gamma_t=\frac{1}{\sqrt{T}}$ and $\bar{\phi}(\hat \x_{0}) - 2f(\hat \x_{0}, \hat{S}_{0})=A_0$,  then

\begin{align}\label{eq:finalS}
\sum_{t=1}^{T}\gamma_t(\bar{\phi}(\hat\x_{t}) - 2f(\hat \x_{t}, \hat S_{t})) &\leq\sum_{t=1}^{T}\frac{1}{\sqrt{T}}(\sum_{t=0}^{t-1}\frac{5G^2}{\sqrt{T}}(1-\frac{2}{k})^{t}+(1-\frac{2}{k})^{t}A_0)
\notag\\&\leq
\sum_{t=1}^{T}\frac{1}{\sqrt{T}}(\frac{5kG^2}{2\sqrt{T}}+(1-\frac{2}{k})^{t}A_0)
\\&\leq {k\beta}
\end{align}
where $\beta=\frac{5kG^2}{2}+\frac{kA_0}{2\sqrt{T}}$
and finally for update of $S$ we get:
\begin{equation}\label{eq:finalSegrg}
\frac{\sum_{t=1}^{T}\gamma_t(\bar{\phi}(\hat \x_{t}) - 2f(\hat \x_{t}, \hat S_{t}))}{\sum_{t=1}^{T}\gamma_t} \leq\frac{k\beta}{{\sqrt{T}}}
\end{equation}

Adding up \eqref{eq:finalSegrg} and \eqref{eq:x EgRg} we have:
%%%%%%%%%%%%%%%%%%

  \begin{align}\label{eq:finalegrg1}
     \frac{\sum_{t=1}^{T}\gamma_t\left [0.5\bar{\phi}(\hat \x_{t})-f(\x,\hat{S}_t)\right]}{\sum_{t=1}^{T}\gamma_t}\leq \frac{K_1}{\sqrt{T}}+\frac{k\beta}{\sqrt{T}}\leq \frac{K}{\sqrt{T}}
    \end{align}
 %%%%%%%%%%%%%%% 
from this for every $S$
%%%%%%%%%%%%%%%%%%

  \begin{align}\label{eq:finalegrg}
     \frac{\sum_{t=1}^{T}\gamma_t\left [0.5f(\hat \x_{t},S)-f(\x,\hat{S}_t)\right]}{\sum_{t=1}^{T}\gamma_t}\leq \frac{K}{\sqrt{T}}
    \end{align}
 %%%%%%%%%%%%%%% 

Similar to \eqref{eq:EGG5} to \eqref{eq:EGG6}, \eqref{eq:finalegrg} results $\x_{sol}=\frac{1}{T}\sum_{t=1}^{T} \hat \x_{t}$, to be $(1/2,\epsilon)$-approximate minimax solution.

 \subsection{Maxmin Result}\label{sec:Proof of Theorem maxmin}

In this section, we introduce maxmin convex-submodular problem and discuss how we can exploit the  algorithms described in the previous sections for the maxmin problem. Formally, consider the function $f:\mathbb{R}^d \times 2^V \to \mathbb{R}_{+}$, where $f(\x,.)$ is submodular for every $\x$ and $f(.,S)$ is convex for every $S$. Then, the maxmin convex-submodular problem is an optimization problem where the maximization is over continuous variable and minimization is over a discrete variable as
\begin{equation} \label{eq:main_problem2}
   {\rm{OPT}}_{maxmin} \triangleq \max_{S \in \mathcal{I}}\min_{x \in \mathcal{X}} f(\x,S),
\end{equation}
 
%%%%%%%%%%%%%%%%%%%%%%%%%%%%
Due to hardness of the max-min problem as we stated in Theorem \ref{thm:saddle_neg} and Appendix \ref{sec:proof of negative result1}, we cannot drive the same result for the maxmin problem as we did for minimax problem. In general, finding an approximation solution for problem \eqref{eq:main_problem2} is NP-hard. Our result as stated in theorem~\ref{thm:maxminresult} proves that $\cup_{t=1}^{T} S_t$ is an approximate solution for \eqref{eq:main_problem2} which has a larger cardinality than our cardinality constraint (at most $Tk$ elements). Although, the set $\cup_{t=1}^{T} S_t$ is not feasible solution, our algorithm converges quickly, and we can use the small number of steps to solve such a problem which means even for small $T$ the set $\cup_{t=1}^{T} S_t$ can solve maxmin problem approximately. This result is similar to the bi-criterion solutions for robust submodular maximization studied in \citep{krause2008robust}, where the authors propose an approach that finds a set that violates the cardinality constraint, but it is within logarithmic factor of the constraint.
\begin{theorem} \label{thm:maxminresult}
Consider all algorithms stated in Algorithms section, if the functions $f$ is convex monotone submodular, and  Assumption~\ref{assum:smooth} holds (and Assumption~\ref{assum:bounded_grad} holds for Gradient Greedy(GG), and Gradient Replacement-greedy(GRG)), then the set $\cup_{t=1}^{T} S_t$   is $({\alpha},\epsilon)-$approximate solution for maxmin convex-submodular problem with cardinality constraint after $\mathcal{O}({1}/{\epsilon^2})$ iterations. Note that parameter $\alpha$ is $\alpha={(1-{1}/{e})}^{-1} $ for Gradient and Extra-gradient Greedy, $\alpha=2 $ for Gradient Replacement-greedy,  $\alpha={2+\frac{k}{k-1}} $ for Extra-gradient Replacement-greedy.
\end{theorem}

%%%%%%%%%%%%%%%%%%%%%%%%%%%%%%%%%%%%%%%%%%

 \subsubsection{Proof of Theorem~\ref{thm:maxminresult} for Gradient Greedy(GG)}
  %%%%%%%%%%%%%%%%%%
 If we let $S^*=\argmax_S\min_\x f(\x,S)$ we know that for every $t$ we have $f(\x_{t-1},S^*)\geq \min_\x f(\x,S^*)=\max_S\min_\x f(\x,S)$. Therefore, if we let $S=S^*$ in \eqref{eq:gg8} we have:
 
  \begin{align}\label{eq:gg12}
({1-\frac{1}{e}})\max_S\min_\x f(\x,S)-\frac{\sum_{t=1}^{T} \alpha (f(\tilde{\x},S_t))}{ \sum_{t=1}^{T}\alpha}\leq \frac{2M^{2}\alpha T+\frac{H}{2\alpha}}{T}
 \end{align}
 
Also, if we let $\hat \x=\argmin_\x f(\x,\cup_t S_t)$ and put $\tilde{\x}=\hat \x$ in \eqref{eq:gg12} then because $f(\hat{\x},S_t)\leq f(\hat{\x},\cup_t S_t)$ we have:
 
   \begin{align}\label{eq:gg13}
({1-\frac{1}{e}})\max_S\min_\x f(\x,S)-\min_\x f(\x,\cup_t S_t)\leq \frac{2M^{2}\alpha T+\frac{H}{2\alpha}}{T}
 \end{align}
 
 and by using $\alpha=\frac{1}{\sqrt{T}}$:
 
 \begin{align}\label{eq:gg14}
({1-\frac{1}{e}})\max_S\min_\x f(\x,S)-\min_\x f(\x,\cup_t S_t)\leq \frac{2M^{2}\alpha T+\frac{H}{2\alpha}}{\sqrt{T}}
 \end{align}
Now, using specific choices $K=2M^{2}+\frac{H}{2}$ and let $T=\frac{K^2}{\epsilon^2}$; we obtain that  $ S_{sol}=\cup_t S_t$ is a $((1-1/e)^{-1},\epsilon)$-approximate maxmin solution.
 %%%%%%%%%%%%%%%%%%%%%%%%%%%%%%%%%%%%%%%%%%%
 \subsubsection{Proof of Theorem~\ref{thm:maxminresult} for Gradient Replacement-greedy(GRG)}
  %%%%%%%%%%%%%%%%%%
  If we let $S^*=\argmax_S\min_\x f(\x,S)$ we know that for every $t$ we have $f(\x_{t-1},S^*)\geq min_\x f(\x,S^*)=\max_S\min_\x f(\x,S)$. Therefore, if we let $S=S^*$ in \eqref{eq:GRG4} we have :
  \begin{align}\label{eq:GRG7}
\frac{1}{2}\max_S\min_\x f(\x,S)-   \frac{1}{T} \sum_{t=1}^{T} f(\tilde{\x},S_t)\leq \frac{K}{\sqrt{T}}
 \end{align}
  Let $\hat \x=\argmin_{\x}f(\x,\cup_t S_t)$ and put $\tilde{\x}=\hat \x$ in \eqref{eq:GRG7} then because $f(\hat{\x},S_t)\leq f(\hat{\x},\cup_t S_t)$ we have:
  \begin{align}\label{eq:GRG71}
\frac{1}{2}\max_S\min_\x f(\x,S)-  \min_{\x}f(\x,\cup_t S_t)&\leq\frac{1}{2}\max_S\min_\x f(\x,S)-   \frac{1}{T} \sum_{t=1}^{T} f(\hat{\x},S_t)\notag\\&\leq \frac{K}{\sqrt{T}}
 \end{align}
 let $T=\frac{K^2}{\epsilon^2}$; then  $ S_{sol}=\cup_t S_t$ is a $(2,\epsilon)$-approximate maxmin solution.
%%%%%%%%%%%%%%%%%%%%%%%%%%%%%%%%%%%%%%%%%%%%%
 \subsubsection{Proof of Theorem~\ref{thm:maxminresult} for Extra-gradient Greedy(EGG)}
  %%%%%%%%%%%%%%%%%%
  If we let $S^*=\argmax_S\min_\x f(\x,S)$ we know that for every $t$ we have $f(\hat \x_{t-1},S^*)\geq min_\x f(\x,S^*)=\max_S\min_\x f(\x,S)$. Therefore, if we let $S=S^*$ in \eqref{eq:EGGmain} we have :
  \begin{align}\label{eq:EGG7}
(1-1/e)\max_S\min_\x f(\x,S)-   \frac{1}{T} \sum_{t=1}^{T} f(\tilde{\x},S_t)\leq \frac{K}{\sqrt{T}}
 \end{align}
  Let $\hat \x=\argmin_{\x} f(\x,\cup_t S_t)$ and put $\tilde{\x}=\hat \x$ in \eqref{eq:GRG7} then because $f(\hat{\x},S_t)\leq f(\hat{\x},\cup_t S_t)$ we have:
  \begin{align}\label{eq:EGG8}
(1-1/e)\max_S\min_\x f(\x,S)&-  \min_\x f(\x,\cup_t S_t)\\&\leq(1-1/e)\max_S\min_\x f(\x,S)-   \frac{1}{T} \sum_{t=1}^{T} f(\hat{\x},S_t)\notag\\&\leq \frac{K}{\sqrt{T}}
 \end{align}
 Let $T=\frac{H^2}{\epsilon^2}$; then, $S_{sol}=\cup_t S_t$ is $((1-1/e)^{-1},\epsilon)$ approximate maxmin solution.
%%%%%%%%%%%%%%%%%%%%%%%%%%%%%%%%%%%%%%%%%%%%%
 \subsubsection{Proof of Theorem~\ref{thm:maxminresult} for Extra-gradient Replacement-greedy (EGRG)}
Similar to \eqref{eq:EGG7}, and \eqref{eq:EGG8}, \eqref{eq:finalegrg} results $S_{sol}=\cup_t S_t$ to be  $((2+\frac{k}{k-1}),\epsilon)$-approximate maxmin solution.

\subsection{Proof of Theorem~\ref{thm:EGcont2}: Extra Gradient on Continuous Extension Convergence}\label{apendix:Extra Gradient on continuous Extension}
In this section, we will focus on convergence analysis of Extra Gradient on continuous extension. We first provide two propositions and matroid definition that will help us in the proof.
{\defi Let $\I$ be a nonempty family of allowable subsets of the ground set $V$, then the tuple $(V, \I)$ is a \textit{matroid} if and only if the following conditions hold:
\begin{enumerate}
    \item For any $A\subset B\subset V$, if $B\in\I$, then $A\in\I$  
    \item For all $A, B\in\I$, if $|A|<|B|$, then there is an $e\in B\backslash A$ such that $A \cup\{e\} \in\I$.
\end{enumerate}
}

\begin{proposition}\label{prop:1}
we have that
\begin{equation}
    {\rm{OPT}} \triangleq \min_{\x \in \mathcal{C}}\max_{S \in \mathcal{I}} f(\x,S) = \min_{\x \in \mathcal{C}} \max_{\y \in \mathcal{K}} F(\x,\y).
\end{equation}

\end{proposition}

Furthermore, the function $F$ has the following properties (assuming differentibility):
\begin{proposition}\label{prop:2}
we have for function $F$(\citep{hassani2017gradient}):
\begin{align*}
    & \forall \x_1, \x_2 \in \mathbb{R}^d:
    F(\x_1,\y) - F(\x_2,\y) \leq \langle \nabla_\x F(\x_1,\y), \x_1 - \x_2 \rangle,
    \\
& \forall \y_1, \y_2 \in \mathbb{R}^d:
    F(\x,\y_2) - 2 F(\x,\y_1) \leq \langle \nabla_\y F(\x,\y_1), \y_2 - \y_1 \rangle.
\end{align*}
\end{proposition}

 using same procedure as Extra-gradient Greedy we drive following equations similar to \eqref{eq:extragradient main}:
%%%%%%%%%%%%%%%%%%%%%%%%%%%%
    \begin{align*}
    \langle-\gamma_t\nabla_\y F(\x_t,\y_t),\hat{\y}_t-\y\rangle\leq \|\y-{\y}_t\|^2-\|\y-\hat{\y}_t\|^2-\|\hat{\y}_t-\y_t\|^2
    \\
    \langle\gamma_t\nabla_\x F(\x_t,\y_t),\hat{\x}_t-\x\rangle\leq \|\x-{\x}_t\|^2-\|\x-\hat{\x}_t\|^2-\|\hat{\x}_t-\x_t\|^2
    \\
    \langle-\gamma_t\nabla_\y F(\hat \x_t,\hat \y_t),{\y}_{t+1}-\y\rangle\leq \|\y-\hat{\y}_t\|^2-\|\y-{\y}_{t+1}\|^2-\|{\y}_{t+1}-\hat \y_t\|^2
   \\
    \langle\gamma_t\nabla_\x F(\hat \x_t,\hat \y_t),{\x}_{t+1}-\x\rangle\leq \|\x-\hat{\x}_t\|^2-\|\x-{\x}_{t+1}\|^2-\|{\x}_{t+1}-\hat \x_t\|^2
    \\
    \langle-\gamma_t\nabla_\y F(\hat \x_t,\hat \y_t),\hat{\y}_{t}-\y_{t+1}\rangle\leq \|\y_{t+1}-{\y}_t\|^2-\|\y-{\y}_{t+1}\|^2-\|{\y}_{t+1}-\y_t\|^2
    \\
    \langle\gamma_t\nabla_\x F(\hat \x_t,\hat \y_t),\hat{\x}_{t}-\x_{t+1}\rangle\leq \|\x_{t+1}-{\x}_t\|^2-\|\x-{\x}_{t+1}\|^2-\|{\x}_{t+1}-\x_t\|^2
    \end{align*}
    %%%%%%%%%%%%%%%%%%%%%%%%%%%%
    combing the above equations we have:
     %%%%%%%%%%%%%%%%%%%%%%%%%%%%  
     
     \begin{align*}
    \langle-\gamma_t\nabla_\y F(\hat \x_t,\hat \y_t),\hat{\y}_{t}-\y\rangle&=  \langle-\gamma_t\nabla_\y F(\hat \x_t,\hat \y_t),\hat{\y}_{t}-\y_{t+1}\rangle\\&\qquad+ \langle-\gamma_t\nabla_\y F(\hat \x_t,\hat \y_t),{\y}_{t+1}-\y\rangle
    \\&=\langle-\gamma_t\nabla_\y F(\hat \x_t,\hat \y_t)+\gamma_t\nabla_\y F( \x_t, \y_t),\hat{\y}_{t}-\y_{t+1}\rangle\\&\qquad+ \langle-\gamma_t\nabla_\y F( \x_t, \y_t),\hat{\y}_{t}-\y_{t+1}\rangle\\&\qquad+\langle-\gamma_t\nabla_\y F(\hat \x_t,\hat \y_t),{\y}_{t+1}-\y\rangle\\&\leq\langle-\gamma_t\nabla_\y F(\hat \x_t,\hat \y_t)+\gamma_t\nabla_\y F( \x_t, \y_t),\hat{\y}_{t}-\y_{t+1}\rangle\\&\qquad-\|\y-{\y}_{t+1}\|^2-\|{\y}_{t+1}-\y_t\|^2\\&\qquad+\|\y-{\y}_t\|^2-\|\y-{\y}_{t+1}\|^2
    \end{align*}
    
 \begin{align*}
    \langle\gamma_t\nabla_\x F(\hat \x_t,\hat \y_t),\hat{\x}_{t}-\x\rangle&=  \langle\gamma_t\nabla_\x F(\hat \x_t,\hat \y_t),\hat{\x}_{t}-\x_{t+1}\rangle\\&\quad+ \langle\gamma_t\nabla_\x F(\hat \x_t,\hat \y_t),{\x}_{t+1}-\x\rangle
    \\&=\langle\gamma_t\nabla_\x F(\hat \x_t,\hat \y_t)-\gamma_t\nabla_\x F( \x_t, \y_t),\hat{\x}_{t}-\x_{t+1}\rangle\\&\quad+ \langle\gamma_t\nabla_\x F( \x_t, \y_t),\hat{\x}_{t}-\x_{t+1}\rangle\\&\quad+\langle\gamma_t\nabla_\x F(\hat \x_t,\hat \y_t),{\x}_{t+1}-\x\rangle\\&\leq\langle\gamma_t\nabla_\x F(\hat \x_t,\hat \y_t)-\gamma_t\nabla_\x F( \x_t, \y_t),\hat{\x}_{t}-\x_{t+1}\rangle\\&\quad-\|\x-{\x}_{t+1}\|^2-\|{\x}_{t+1}-\x_t\|^2\\&\quad+\|\x-{\x}_t\|^2-\|\x-{\x}_{t+1}\|^2
    \end{align*}
    
    %%%%%%%%%%%%%%%%%%%%%%%%%%%%
    let 
    \begin{align*}
        \sigma_t^\x =& \langle\gamma_t\nabla_\x F(\hat \x_t,\hat \y_t)-\gamma_t\nabla_\x F( \x_t, \y_t),\hat{\x}_{t}-\x_{t+1}\rangle-\|\x-{\x}_{t+1}\|^2\\&-\|{\x}_{t+1}-\x_t\|^2
    \end{align*}
    and
    \begin{align*}\sigma_t^\y =& \langle-\gamma_t\nabla_\y F(\hat \x_t,\hat \y_t)+\gamma_t\nabla_\y F( \x_t, \y_t),\hat{\y}_{t}-\y_{t+1}\rangle-\|\y-{\y}_{t+1}\|^2\\&-\|{\y}_{t+1}-\y_t\|^2
    \end{align*}

    then $\sigma_t^\x\leq 0$, $\sigma_t^\y\leq 0$
    if $\gamma_t\leq \frac{1}{\max\{L_\x,L_\y\}}$(check \citep{nemirovski2004prox} for more details)
    ; which results in:
    %%%%%%%%%%%%%%%%%%%%%%%%%%%%
     
     \begin{align}\label{eq:mirrorFinal1x}
    \langle-\gamma_t\nabla_\y F(\hat \y_t,\hat \y_t),\hat{\y}_{t}-\y\rangle& \leq \|\y-{\y}_t\|^2-\|\y-{\y}_{t+1}\|^2
    \end{align}
     \begin{align}\label{eq:mirrorFinal2y}
    \langle\gamma_t\nabla_\x F(\hat \x_t,\hat \y_t),\hat{\x}_{t}-\x\rangle& \leq \|\x-{\x}_t\|^2-\|\x-{\x}_{t+1}\|^2
    \end{align}
    %%%%%%%%%%%%%%%%%%%%%%%%%%%%
    combing above equations with proposition \ref{prop:2} we have:
 %%%%%%%%%%%%%%%%%%%%%%%%%%%%
  \begin{align}\label{eq:mirrorFinalx}
    2\gamma_tF(\hat \x_{t},\hat \y_{t})-2\gamma_tF( \x,\hat \y_{t})&\leq2\langle\gamma_t\nabla_\x F(\hat \x_t,\hat \y_t),\hat{\x}_{t}-x\rangle\notag \\&\leq 2\|\x-{\x}_t\|^2-2\|\x-{\x}_{t+1}\|^2
    \end{align}
    %%%%%%%%%%%%%%%%%%%%%%%%%%%%
     \begin{align}\label{eq:mirrorFinaly}
     -2\gamma_tF(\hat \x_{t},\hat \y_{t})+\gamma_tF(\hat \x_{t} ,\y)&\leq\langle-\gamma_t\nabla_\y F(\hat \x_t,\hat \y_t),\hat{\y}_{t}-\y\rangle\notag \\& \leq \|\y-{\y}_t\|^2-\|\y-{\y}_{t+1}\|^2
    \end{align}
    %%%%%%%%%%%%%%%%%%%%%%%%%%%%
     \begin{align}\label{eq:mirrorFinalyx}
     -2\gamma_tF( \x,\hat \y_{t})+\gamma_tF(\hat \x_{t} ,\y)&\leq \|\y-{\y}_t\|^2-\|\y-{\y}_{t+1}\|^2\notag \\&\quad+2\|\x-{\x}_t\|^2-2\|\x-{\x}_{t+1}\|^2
    \end{align}
    %%%%%%%%%%%%%%%%%%%%%%%%%%%%
    summing over $t$ in \eqref{eq:mirrorFinalyx} and divide both side by $\sum_{t=1}^{T}\gamma_t $ (set of variable $\x$ and $\y$ is bounded i.e. $\|\y\|^2\leq H,\|\x\|^2\leq H$):
     %%%%%%%%%%%%%%%%%%%%%%%%%%%%
      \begin{align}
     \frac{\sum_{t=1}^{T}\gamma_t\left [-2F( \x,\hat \y_{t})+F(\hat \x_{t} ,\y)\right]}{\sum_{t=1}^{T}\gamma_t}\leq\frac{\|\y-{\y}_1\|^2+2\|\x-{\x}_1\|^2}{\sum_{t=1}^{T}\gamma_t}\leq \frac{3H}{\gamma T}
    \end{align}
which means same as before let $T=\frac{\sqrt{3H}}{\gamma\epsilon}$ and constant step size $\gamma_t=
\gamma$, and $\x^{*}=\argmin\max f(\x,\y)$ we have:

  \begin{align}
 \frac{1}{2}f(  \frac{1}{T} \sum_{t=1}^{T} \hat \x_t, \y)-  \min_\x\max_\y F(\x,\y)\leq \frac{1}{2}F(  \frac{1}{T} \sum_{t=1}^{T}  \x_t, \y)-  \frac{1}{T} \sum_{t=1}^{T} f({\x}^{*},\y^t)\leq 
 \epsilon
 \end{align}
then using proposition~\ref{prop:1}, $\x_{sol}=\frac{1}{T}\sum_{t=1}^{T} \hat \x_{t}$ is $(0.5,\epsilon)$-approximate minimax solution.

\end{document}